\documentclass[10pt]{article}%
\usepackage{amsmath}
\usepackage{amsfonts}
\usepackage{amssymb}
\usepackage{graphicx}%
\usepackage{mathrsfs}
\usepackage{xcolor}
\usepackage[body={15.5cm,21cm}, top=3cm]{geometry}
\DeclareMathOperator*{\essinf}{ess\,inf}
\providecommand{\U}[1]{\protect \rule{.1in}{.1in}}
\newtheorem{theorem}{Theorem}[section]

\newtheorem{definition}[theorem]{Definition}

\newtheorem{remark}[theorem]{Remark}

\newenvironment{proof}[1][Proof]{\noindent \textbf{#1.} }{\  \rule{0.5em}{0.5em}}
\allowdisplaybreaks[2]
\begin{document}

\title{Diagonally quadratic BSDE with oblique reflection and optimal switching}
\author{Peng Luo \thanks{School of Mathematical Sciences, Shanghai Jiao Tong University, China (peng.luo@sjtu.edu.cn). Financial support from the National Natural Science Foundation of China (Grant No. 12101400) is gratefully acknowledged.}
\and Mengbo Zhu \thanks{School of Mathematical Sciences, Ocean University of China, China (21201131039@stu.ouc.edu.cn)}}

\maketitle

\begin{abstract}
The present paper is devoted to the study of diagonally quadratic backward stochastic differential equation with oblique reflection. Using a penalization approach, we show the existence fo a solution by providing some delicated a priori estimates. We further obtain the uniqueness by verifying the first component of the solution is indeed the value of a switching probelm for quadratic BSDEs. Moreover, we provide an extension for the solvability and apply our results to study a risk-sensitive switching problem for functional stochastic differential equations. 
\end{abstract}


\textbf{MSC-classification}: 60H10, 93H20.
\section{Introduction}
Given a multidimensional Brownian motion on a complete probability space $(\Omega,\mathcal{F},P)$, we consider the following reflected backward stochastic differential equation (RBSDE for short) with oblique reflection
\begin{equation}\label{eq:intro_BSDE}
    \begin{cases}
    &Y_i(t)=\xi_i+\int_t^Tg_i(s,Y_i(s),Z_i(s))ds-\int_t^TdK_i(s)-\int_t^TZ_i(s)dW(s),~~t\in[0,T],\\
    &Y_i(t)\leq \min\limits_{j\neq i}\{Y_j(t)+k_{i,j}\},\\
    &\int_0^T\left(Y_i(s)-\min\limits_{j\neq i}\left\{Y_j(s)+k_{i,j}\right\}\right)dK_i(s)=0, ~~i=1,2,\ldots,n,
    \end{cases}
\end{equation}
where $\xi$ is an $n$-dimensional measurable random variable, $T>0$ is a finite time horizon, $g$ is a given function and $k_{i,j}$ are given constants for all $1\leq i,j\leq n$. In this paper, we provide conditions under which RBSDE \eqref{eq:intro_BSDE} has a unique solution in the case where $g$ is diagonally quadratic in $Z$. 

When $g$ is uniformly in $Y$ and $Z$, solvability of RBSDE \eqref{eq:intro_BSDE} has been studied in Hu and Tang \cite{HT1}. They also provided an application in studing a switching problem for stochastic differential equations of functional type and a probabilistic representation for the viscosity solution of a system of variational inequalities.  By assuming a monotonicity condition, Hamad\`{e}ne and Zhang \cite{HZ} condsidered RBSDEs with more general oblique constraints and generators depending on the whole vector $Y$. In consideration of linear reflections, Chassagneux et al. \cite{CE} further obtained the existence and uniquess of the solution without the monotonicity condition. Chassagneux and Richou \cite{CR} futher studied existence and uniqueness of solutions for multidimensional RBSDEs in a non-empty open convex domain. More recently, B\'{e}n\'{e}zet et al. \cite{BCR} considered a new class of obliquely reflected BSDEs which is related to optimal switching problems with controlled randomisation. On the other hand, obliquely reflected BSDEs are closely related to optimal switching problems which arise naturally in energy pricing and investment timing modles, see \cite{CL,Lud} and references therein. In particular, Porchet et al. \cite{PTW} studied a system of Markovian RBSDE where the generator $g$ is quadratic in $Z$.

This work is mainly motivated by the recent developments of multidimensional quadratic BSDEs (see \cite{Te,HT,JKL,Luo1,CN,XZ}). We consider diagonally quadratic BSDE with oblique reflection. We first introduce a sequence of penalized BSDEs and provide delicated a priori estimates for their solutions. Thanks to the comparion result for diagonally quadratic BSDEs recently obtained in \cite{Luo}, we show that the convergence of the value processes of these penalized BSDEs by a monotone convergence argument. Based on the a priori estimates, we are able to show that the penalized solutions indeed converge in some suitable spaces by heavily using the properties of BMO martingales, which yields the existence of a solution. We further verify that this solution is in some finer space. Similar to Hu and Tang \cite{HT1}, we obtain the uniqueness by showing that the first component of the solution is the value of a switching problem for quadratic BSDEs. Following a similar argument as in \cite{CE}, we extend the existence and uniqueness results to a general case where the generator can depend on the whole vector of the value process. As an application, we study a risk-sensitive switching problem of functional SDEs which is solved by using the obtained solvability result for diagonally quadratic BSDE with oblique reflection.

The remainder of the paper is organized as follows. In the next section, we present our settings and main results. In subsection 2.1, we give the existence result while the uniqueness result is stated in subsection 2.2. Subsection 2.3 is devoted to the proof of Theorem \ref{thm:main-ex}. The risk-sensitive switching problem of functional SDEs is studied in Section 3.

\section{Diagonally quadratic BSDE with oblique reflection}
Let $W=(W_t)_{t\geq 0}$ be a $d$-dimensional Brownian motion on a probability space $(\Omega, {\cal F}, P)$. Let $(\mathcal{F}_t)_{t\geq 0}$ be the augmented filtration generated by $W$. Throughout, we fix a $T\in (0,\infty)$. We endow $\Omega \times [0,T]$ with the predictable $\sigma$-algebra $\mathcal{P}$ and $\mathbb{R}^n$ with its Borel $\sigma$-algebra $\mathcal{B}(\mathbb{R}^n)$. Equalities and inequalities between random variables and processes are understood in the $P$-a.s. and $P\otimes dt$-a.e. sense, respectively. The Euclidean norm is denoted by $|\cdot|$ and $\|\cdot\|_\infty$ denotes the $L^\infty$-norm. $\mathcal{C}_{T}(\mathbb{R}^n)$ denotes the set $C([0,T];\mathbb{R}^n)$ of continuous functions from $[0,T]$ to $\mathbb{R}^n$. For $p>1$, we denote by
\begin{itemize}
\item $\mathcal{S}^p(\mathbb{R}^n)$ (resp. $\mathcal{S}^p_c(\mathbb{R}^n)$) the set of $n$-dimensional c\`{a}dl\`{a}g (resp. continuous) adapted processes $Y$ on $[0,T]$ such that
\begin{equation*}
\|Y\|_{\mathcal{S}^p}:=E\left[\sup_{0\leq t\leq T} |Y(t)|^p\right]^{\frac{1}{p}}< \infty;
\end{equation*}
\item $\mathcal{S}^\infty(\mathbb{R}^n)$ (resp. $\mathcal{S}^\infty_c(\mathbb{R}^n)$) the set of $n$-dimensional c\`{a}dl\`{a}g (resp. continuous) adapted processes $Y$ on $[0,T]$ such that
\begin{equation*}
\|Y\|_{\mathcal{S}^\infty}:=\bigg\| \sup_{0\leq t\leq T} |Y(t)| \bigg\|_{\infty} < \infty;
\end{equation*}
\item  $\mathcal{A}^p(\mathbb{R}^n)$ (resp.  $\mathcal{A}^p_c(\mathbb{R}^n)$) is the closed subset of $\mathcal{S}^{p}(\mathbb{R}^n)$ (resp. $\mathcal{S}^{p}_c(\mathbb{R}^n)$) consisting of processes $K = (K(t))_{0\leq t\leq T}$ starting from the origin such that $K_i(t)$ is non-decreasing in $t$ for all $1\leq i\leq n$;
\item $\mathcal{A}^{\infty}(\mathbb{R}^n)$ (resp.  $\mathcal{A}^{\infty}_c(\mathbb{R}^n)$) is the closed subset of $\mathcal{S}^{\infty}(\mathbb{R}^n)$ (resp. $\mathcal{S}^{\infty}_c(\mathbb{R}^n)$) consisting of processes $K = (K(t))_{0\leq t\leq T}$ starting from the origin such that $K_i(t)$ is non-decreasing in $t$ for all $1\leq i\leq n$;
\item $\mathcal{H}^p(\mathbb{R}^{n\times d})$ the set of predictable $\mathbb{R}^{n\times d}$-valued processes $Z$ such that
    \begin{equation*}
    \|Z\|_{\mathcal{H}^p}=E\left[\left(\int_0^T|Z(s)|^2ds\right)^{\frac{p}{2}}\right]^{\frac{1}{p}}<\infty.
    \end{equation*}
\end{itemize}
Let $\mathcal{T}$ be the set of all stopping times with values in $[0,T]$. For any uniformly integrable martingale $M$ with $M_0=0$ and for $p\geq 1$, we set
\begin{equation*}
\|M\|_{BMO_p(P)}:=\sup_{\tau\in\mathcal{T}}\|E[|M(T)-M(\tau)|^p|\mathcal{F}_\tau]^{\frac{1}{p}}\|_\infty.
\end{equation*}
The class $\{M:~\|M\|_{BMO_p}<\infty\}$ is denoted by $BMO_p$, which is written as $BMO(P)$ when it is necessary to indicate the underlying probability measure $P$. In particular, we will denote it by $BMO$ when $p=2$. For $(\alpha\cdot W)_t:=\int_0^t\alpha(s)dW(s)$ in $BMO$, the corresponding stochastic exponential is denoted by $\mathcal{E}_{t}(\alpha\cdot W)$.

We consider the following reflected backward stochastic differential equation with oblique reflection: 
\begin{equation}\label{BSDE}
\begin{cases}
&Y_i(t)=\xi_i+\int_t^Tg_i(s,Y_i(s),Z_i(s))ds-\int_t^TdK_i(s)-\int_t^TZ_i(s)dW(s),~~t\in[0,T],\\
&Y_i(t)\leq \min\limits_{j\neq i}\{Y_j(t)+k_{i,j}\},\\
&\int_0^T\left(Y_i(s)-\min\limits_{j\neq i}\left\{Y_j(s)+k_{i,j}\right\}\right)dK_i(s)=0, ~~i=1,2,\ldots,n,
\end{cases}
\end{equation}
where $\xi$ is $\mathbb{R}^n$-valued and $\mathcal{F}_T$-measurable random variable which is bounded, $k_{i,j}\in\mathbb{R}$ for all $i,j=1,2,\cdots,n$. Let $\gamma>0$ be a constant, we will make the following assumptions:
\begin{itemize}
   \item[(A1)] For any $i=1,2,\cdots,n$, $g_i:\Omega\times[0,T]\times\mathbb{R}\times\mathbb{R}^d\rightarrow\mathbb{R}$ satisfies that $g_i(\cdot,y,z)$ is adapted for each $y\in\mathbb{R},z\in\mathbb{R}^d$. Moreover, it holds that $g_i(t,y,z)=h_i(t,y)+f(t,z)$, where $h_i$ and $f$ satisfy
       \begin{align*}
        &|h_i(t,y)|\leq \gamma(1+|y|),~~|h_i(t,y)-h_i(t,\bar{y})|\leq \gamma|y-\bar{y}|,~~~y,\bar{y}\in\mathbb{R};\\
       &|f(t,z)|\leq \gamma(1+|z|^2),~~|f(t,z)-f(t,\bar{z})|\leq \gamma(1+|z|+|\bar{z}|)|z-\bar{z}|,~~~z,\bar{z}\in\mathbb{R}^d;
       \end{align*}
   \item[(A2)]
       For any $i,j=1,2,\cdots,n$, $k_{i,j}$ satisfies that $k_{i,i}=0$ and $k_{i,j}>0$ if $i\neq j$;
   \item[(A3)] For any $i,j,l=1,2,\cdots,n$, it holds that $k_{i,j}+k_{j,l}\geq k_{i,l}$;
   \item[(A4)] For any $i,j,l=1,2,\cdots,n$, it holds that $k_{i,j}+k_{j,l}> k_{i,l}$.
   \end{itemize}
   We also introduce the following domain $Q$:
   \begin{equation*}
      Q:=\left\{y\in\mathbb{R}^n:~y_i<y_j+k_{i,j}~\text{for any}~i,j=1,2,\cdots,n~\text{satisfying}~j\neq i\right\}
   \end{equation*}
   which is convex and bounded. The closure of $Q$ is denoted by $\overline{Q}$.

   Our first main result gives the existence and uniquness of the solution for RBSDE \eqref{BSDE} which follows from Theorem \ref{thm:exist} and Theorem \ref{thm:uniqueness}.
   \begin{theorem}\label{thm:main}
    Assume that (A1),(A2),(A4) hold and $\xi$ takes values in $\overline{Q}$. Then RBSDE \eqref{BSDE} admits a unique solution $(Y,Z,K)$ such that $(Y,Z\cdot W,K)\in\mathcal{S}^{\infty}_c(\mathbb{R}^n)\times BMO \times\mathcal{A}^{\infty}_c(\mathbb{R}^n)$.
   \end{theorem}
   \begin{remark}
       In oder to obtain our solvability result, we require the generator to satisfy $g_i(t,y,z)=h_i(t,y)+f(t,z)$, which is also needed for the converse comparison theorem for diagonally quadratic BSDEs obtained in \cite{Luo}. 
   \end{remark}
   \begin{remark}
       The existence part holds under weaker assumption than the uniqueness part. Actually, under the assumptions (A1),(A2) and (A3), RBSDE \eqref{BSDE} has a solution, see Theorem \ref{thm:exist}.
   \end{remark}
   Next, we would like to extend Theorem \ref{thm:main} to a general case where the generator is allowed to depend on the whole vector of the value process. Indeed, we consider the following RBSDE
   \begin{equation}\label{BSDE1}
    \begin{cases}
    &Y_i(t)=\xi_i+\int_t^Tg_i(s,Y(s),Z_i(s))ds-\int_t^TdK_i(s)-\int_t^TZ_i(s)dW(s),~~t\in[0,T],\\
    &Y_i(t)\leq \min\limits_{j\neq i}\{Y_j(t)+k_{i,j}\},\\
    &\int_0^T\left(Y_i(s)-\min\limits_{j\neq i}\left\{Y_j(s)+k_{i,j}\right\}\right)dK_i(s)=0,~~~i=1,2,\ldots,n,
    \end{cases}
    \end{equation}
    We further introduce the following condition.
   \begin{itemize}
       \item[(A5)] For any $i=1,2,\cdots,n$, $g_i:\Omega\times[0,T]\times\mathbb{R}^n\times\mathbb{R}^d\rightarrow\mathbb{R}$ satisfies that $g_i(\cdot,y,z)$ is adapted for each $y\in\mathbb{R}^n,z\in\mathbb{R}^d$. Moreover, it holds that $g_i(t,y,z)=h_i(t,y)+f(t,z)$, where $h_i$ and $f$ satisfy
       \begin{align*}
        &|h_i(t,y)|\leq \gamma(1+|y|),~~|h_i(t,y)-h_i(t,\bar{y})|\leq \gamma|y-\bar{y}|,~~~y,\bar{y}\in\mathbb{R}^n;\\
       &|f(t,z)|\leq \gamma(1+|z|^2),~~|f(t,z)-f(t,\bar{z})|\leq \gamma(1+|z|+|\bar{z}|)|z-\bar{z}|,~~~z,\bar{z}\in\mathbb{R}^d.
       \end{align*}
   \end{itemize}
   Following a similar technique as in \cite{CE}, we obtain the solvability under this general setting. The proof is given in subsection \ref{extension}.
   \begin{theorem}\label{thm:main-ex}
    Assume that (A2),(A4),(A5) hold and $\xi$ takes values in $\overline{Q}$. Then RBSDE \eqref{BSDE1} admits a unique solution $(Y,Z,K)$ such that $(Y,Z\cdot W,K)\in\mathcal{S}^{\infty}_c(\mathbb{R}^n)\times BMO \times\mathcal{A}^{\infty}_c(\mathbb{R}^n)$.
   \end{theorem}
   \begin{remark}
       Our results still hold true if the constants $k_{i,j}$ are replaced by continuous processes $(k_{i,j}(t))_{0\leq t\leq T}$ by using a similar argument as in the proof of Proposition 3.1 of \cite{CE}.
   \end{remark} 
\subsection{Existence}
The subsection is devoted to show the existence part of Theorem \ref{thm:main}. We state the existence result in the following theorem. 
\begin{theorem}\label{thm:exist}
Assume that (A1),(A2),(A3) hold and $\xi$ takes values in $\overline{Q}$. Then RBSDE \eqref{BSDE} has a solution $(Y,Z,K)$ such that $(Y,Z\cdot W,K)\in\mathcal{S}^{\infty}_c(\mathbb{R}^n)\times BMO \times\mathcal{A}^{\infty}_c(\mathbb{R}^n)$.
\end{theorem}
\begin{proof} The proof is divided into several steps. 

\noindent \textbf{Step 1.} For any $i=1,2,\cdots,n$ and nonnegative integer $m$, we define function $\varphi^{m}_i$ as follows:
\begin{equation*}
\varphi^{m}_i(t,y,z):=g_{i}(t,y_i,z_i)-m\sum_{l=1}^n\left(y_i-y_l-k_{i,l}\right)^+,~~0\leq t\leq T,~y\in\mathbb{R}^n,~z\in\mathbb{R}^{n\times d}.
\end{equation*}
 We introduce the following penalized BSDE:
\begin{equation}\label{BSDE_penalized}
Y^{m}_i(t)=\xi_i+\int_t^T\varphi^m_i(s,Y^{m}(s),Z^{m}_i(s))ds-\int_t^TZ^{m}_i(s)dW(s),~t\in[0,T],~i=1,2,\ldots,n.
\end{equation}
Since (A1) holds, it follows from \cite[Theorem 2.3]{HT} that the penalized BSDE \eqref{BSDE_penalized} admits a unique solution $\left(Y^m,Z^{m}\cdot W\right)\in\mathcal{S}^{\infty}_c\left(\mathbb{R}^n\right)\times BMO$.
On the other hand, for any $i=1,2,\ldots,n$, $t\in[0,T]$, $y\in\mathbb{R}^n$ and $z_i\in\mathbb{R}^d$, it holds that 
\begin{equation*}
\varphi^{m}_i(t,y+\tilde{y},z_i)\geq\varphi^{m+1}_i(t,y,z_i)
\end{equation*}
for any $\tilde{y}\in\mathbb{R}^n$ such that each component $\tilde{y}_j\geq 0$ and the $i$-th component $\tilde{y}_i=0$. Therefore, it follows from \cite[Theorem 2.1]{Luo} that for each $m$ and $t\in[0,T]$,
\begin{equation*}
Y^{m}_i(t)\geq Y^{m+1}_i(t),~~~i=1,2,\ldots,n.
\end{equation*}

Defining
\begin{equation*}
\tilde{\xi}=-\sum_{j=1}^n|\xi_j|~\text{and}~\tilde{g}(t,y,z)=-\sum_{j=1}^n|g_j(t,y,z)|,~~~t\in[0,T],y\in\mathbb{R},z\in\mathbb{R}^d
\end{equation*}
and noting that (A1) holds, it follows from \cite[Theorem 2.2]{BE1} that the following BSDE
\begin{equation*}
\tilde{Y}(t)=\tilde{\xi}+\int_t^T\tilde{g}(s,\tilde{Y}(s),\tilde{Z}(s))ds-\int_t^T\tilde{Z}(s)dW(s),
\end{equation*}
has a unique solution $(\tilde{Y},\tilde{Z}\cdot W)\in\mathcal{S}^{\infty}_c(\mathbb{R})\times BMO$. For any $i=1,\ldots,n$ and $t\in[0,T]$, we define
\begin{equation*}
\bar{Y}_i(t)=\tilde{Y}(t),~~~\bar{Z}_i(t)=\tilde{Z}(t).
\end{equation*}
 Since (A2) holds, obviously we get $\left(\bar{Y}^i-\bar{Y}^l-k_{i,l}\right)^+=0$. Therefore for any $i=1,2,\ldots,n$ and all $m$, $(\bar{Y}_i,\bar{Z}_i)$ satisfies
\begin{equation*}
\bar{Y}_i(t)=\tilde{\xi}+\int_t^T\left(\tilde{g}(s,\bar{Y}_i(s),\bar{Z}_i(s))-m\sum_{l=1}^n\left(\bar{Y}_i(s)-\bar{Y}_l(s)-k_{i,l}\right)^{+}\right)ds-\int_t^T\bar{Z}_i(s)dW(s),~~~~t\in[0,T].
\end{equation*}
Once again, it follows from \cite[Theorem 2.1]{Luo} that for each $m$ and $t\in[0,T]$,
\begin{equation*}
Y^{m}_i(t)\geq \bar{Y}_{i}(t)=\tilde{Y}(t).
\end{equation*}
Hence there exists a constant $C_1>0$ which is independent of $m$ such that
\begin{equation}\label{eq:Yuniform_bound}
    \|Y^m\|_{\mathcal{S}^{\infty}}\leq C_1.
    \end{equation}
For any $i=1,2,\ldots,n$ and all $m$, applying It\^{o} formula to $e^{2(\gamma+1) Y^{m}_i(t)}$ implies that for all $t\in[0,T]$
\begin{align*}
    e^{2(\gamma+1)Y^{m}_i(t)}&=e^{2(\gamma+1)\xi_i}-\int_{t}^T2(\gamma+1)e^{2(\gamma+1)Y^m_i(s)}Z^m_i(s)dW(s)\\
    &+\int_t^T\left(2(\gamma+1)e^{2(\gamma+1)Y^{m}_i(s)}\varphi^m_i\left(s,Y^m(s),Z^m_i(s)\right)-2(\gamma+1)^2e^{2(\gamma+1)Y^m_i(s)}|Z^m_i(s)|^2\right)ds\\
    &\leq e^{2(\gamma+1)\xi_i}-\int_{t}^T2(\gamma+1)e^{2(\gamma+1)Y^m_i(s)}Z^m_i(s)dW(s)\\
    &+\int_t^T\left(2\gamma(\gamma+1)e^{2(\gamma+1)Y^{m}_i(s)}\left(2+|Y^m_i(s)|+|Z^m_i(s)|^2\right)-2(\gamma+1)^2e^{2(\gamma+1)Y^m_i(s)}|Z^m_i(s)|^2\right)ds\\
    &\leq e^{2(\gamma+1)\xi_i}-\int_{t}^T2(\gamma+1)e^{2(\gamma+1)Y^m_i(s)}Z^m_i(s)dW(s)\\
    &+\int_t^T\left(2\gamma(\gamma+1)e^{2(\gamma+1)Y^{m}_i(s)}\left(2+|Y^m_i(s)|\right)-2(\gamma+1)e^{2(\gamma+1)Y^m_i(s)}|Z^m_i(s)|^2\right)ds.
\end{align*}
Hence, it holds that
\begin{align*}
    2(\gamma+1)\int_t^Te^{2(\gamma+1)Y^m_i(s)}|Z^m_i(s)|^2ds&\leq e^{2(\gamma+1)\xi_i}-\int_{t}^T2(\gamma+1)e^{2(\gamma+1)Y^m_i(s)}Z^m_i(s)dW(s)\\
    &+\int_t^T2\gamma(\gamma+1)e^{2(\gamma+1)Y^{m}_i(s)}\left(2+|Y^m_i(s)|\right)ds.
\end{align*}
Taking conditional expectation with respect to $\mathcal{F}_t$ and noting the uniform boundedness of $Y^m$, it holds that
\begin{align*}
    E\left[\int_t^T|Z^m_i(s)|^2ds\bigg|\mathcal{F}_t\right]\leq \frac{1}{2}e^{2(\gamma+1)\|Y^m\|_{\mathcal{S}^{\infty}}}\left(e^{2(\gamma+1)\|\xi\|_{\infty}}+2\gamma(\gamma+1)e^{2(\gamma+1)\|Y^m\|_{\mathcal{S}^{\infty}}}\left(2+\|Y^m\|_{\mathcal{S}^{\infty}}\right)T\right)
\end{align*}
which implies that there exists a constant $C_2>0$ which is independent of $m$ such that
\begin{equation}\label{eq:Zuniform_bound}
\|Z^m\cdot W\|_{BMO}\leq C_2.
\end{equation}
\textbf{Step 2.} We now show that there exists constant $C_3>0$ (independent of $m$) such that
\begin{equation}\label{eq:uniform_penalized}
    \left\|\left(Y^{m}_i-Y^{m}_j-k_{i,j}\right)^{+}\right\|_{\mathcal{S}^{\infty}}\leq\frac{C_3}{m}.
    \end{equation}
Indeed, for any $p>2$, $i,j=1,2,\ldots,n$ and all $m$, applying It\^{o}'s formula to $\left(\left(Y^{m}_i(t)-Y^{m}_j(t)-k_{i,j}\right)^{+}\right)^p$, we have
\begin{align*}
&\left(\left(Y^{m}_i(t)-Y^{m}_j(t)-k_{i,j}\right)^{+}\right)^p+pm\int_t^T\left(\left(Y^{m}_i(s)-Y^{m}_j(s)-k_{i,j}\right)^{+}\right)^pds\\
&\qquad +\frac{p(p-1)}{2}\int_t^T\left(\left(Y^{m}_i(s)-Y^{m}_j(s)-k_{i,j}\right)^{+}\right)^{p-2}|Z^{m}_i(s)-Z^{m}_j(s)|^2ds\\
&=p\int_t^T\left(\left(Y^{m}_i(s)-Y^{m}_j(s)-k_{i,j}\right)^{+}\right)^{p-1}\left(g_i(s,Y^{m}_i(s),Z^{m}_i(s))-g_j(s,Y^{m}_j(s),Z^{m}_j(s))\right)ds\\
&\qquad -p\int_t^T\left(\left(Y^{m}_i(s)-Y^{m}_j(s)-k_{i,j}\right)^{+}\right)^{p-1}\left(Z^{m}_i(s)-Z^{m}_j(s)\right)dW(s)\\
&\qquad
+pm\int_t^T\left(\left(Y^{m}_i(s)-Y^{m}_j(s)-k_{i,j}\right)^{+}\right)^{p-1}\left(Y^{m}_j(s)-Y^{m}_i(s)-k_{j,i}\right)^{+}ds\\
&\qquad
+pm\sum_{j\neq i,l\neq j}\int_t^T\left(\left(Y^{m}_i(s)-Y^{m}_j(s)-k_{i,j}\right)^{+}\right)^{p-1}\left(\left(Y^{m}_j(s)-Y^{m}_l(s)-k_{j,l}\right)^{+}-\left(Y^{m}_i(s)-Y^{m}_l(s)-k_{i,l}\right)^{+}\right)ds.
\end{align*}
For any $i,j=1,2,\ldots,n$, it follows from (A2) that
\begin{equation*}
    k_{i,j}+k_{j,i}\geq 0,
    \end{equation*}
which implies that 
\begin{equation*}
    \left\{y\in\mathbb{R}^n:~y_i-y_j-k_{i,j}>0,~y_j-y_i-k_{j,i}>0\right\}=\emptyset.
    \end{equation*}
Therefore, we have 
\begin{equation*}
\left(\left(Y^{m}_i(s)-Y^{m}_j(s)-k_{i,j}\right)^{+}\right)^{p-1}\left(Y^{m}_j(s)-Y^{m}_i(s)-k_{j,i}\right)^{+}=0,~i,j=1,2,\cdots,n.
\end{equation*}
On the other hand, for any $i,j,l=1,2,\ldots,n$, using (A3)
and the elementary inequality that $x^+-y^+\leq (x-y)^+$ for any two real numbers $x$ and $y$, we have
\begin{align*}
&\left(\left(Y^{m}_i(s)-Y^{m}_j(s)-k_{i,j}\right)^{+}\right)^{p-1}\left(\left(Y^{m}_j(s)-Y^{m}_l(s)-k_{j,l}\right)^{+}-\left(Y^{m}_i(s)-Y^{m}_l(s)-k_{i,l}\right)^{+}\right)\\
&\leq \left(\left(Y^{m}_i(s)-Y^{m}_j(s)-k_{i,j}\right)^{+}\right)^{p-1}\left(Y^{m}_j(s)-Y^{m}_i(s)-k_{j,l}+k_{i,l}\right)^{+}\\
&\leq \left(\left(Y^{m}_i(s)-Y^{m}_j(s)-k_{i,l}+k_{j,l}\right)^{+}\right)^{p-1}\left(Y^{m}_j(s)-Y^{m}_i(s)-k_{j,l}+k_{i,l}\right)^{+}\\
&=\left(\left(Y^{m}_i(s)-Y^{m}_j(s)-k_{i,l}+k_{j,l}\right)^{+}\right)^{p-1}\left(Y^{m}_i(s)-Y^{m}_j(s)-k_{i,l}+k_{j,l}\right)^{-}\\
&=0
\end{align*}
Combing the above, in view of (A1), we have
\begin{align*}
&\left(\left(Y^{m}_i(t)-Y^{m}_j(t)-k_{i,j}\right)^{+}\right)^p+pm\int_t^T\left(\left(Y^{m}_i(s)-Y^{m}_j(s)-k_{i,j}\right)^{+}\right)^pds\\
&\qquad +\frac{p(p-1)}{2}\int_t^T\left(\left(Y^{m}_i(s)-Y^{m}_j(s)-k_{i,j}\right)^{+}\right)^{p-2}|Z^{m}_i(s)-Z^{m}_j(s)|^2ds\\
&\leq p\int_t^T\left(\left(Y^{m}_i(s)-Y^{m}_j(s)-k_{i,j}\right)^{+}\right)^{p-1}\left(g_i(s,Y^{m}_i(s),Z^{m}_i(s))-g_j(s,Y^{m}_j(s),Z^{m}_j(s))\right)ds\\
&\qquad -p\int_t^T\left(\left(Y^{m}_i(s)-Y^{m}_j(s)-k_{i,j}\right)^{+}\right)^{p-1}\left(Z^{m}_i(s)-Z^{m}_j(s)\right)dW_s\\
&= p\int_t^T\left(\left(Y^{m}_i(s)-Y^{m}_j(s)-k_{i,j}\right)^{+}\right)^{p-1}\left(h_i(s,Y^{m}_i(s))-h_j(s,Y^{m}_j(s))\right)ds\\\
&\qquad +p\int_t^T\left(\left(Y^{m}_i(s)-Y^{m}_j(s)-k_{i,j}\right)^{+}\right)^{p-1}\left(f(s,Z^{m}_i(s))-f(s,Z^{m}_j(s))\right)ds\\
&\qquad -p\int_t^T\left(\left(Y^{m}_i(s)-Y^{m}_j(s)-k_{i,j}\right)^{+}\right)^{p-1}\left(Z^{m}_i(s)-Z^{m}_j(s)\right)dW(s)\\
&= p\int_t^T\left(\left(Y^{m}_i(s)-Y^{m}_j(s)-k_{i,j}\right)^{+}\right)^{p-1}\left(h_i(s,Y^{m}_i(s))-h_j(s,Y^{m}_j(s))\right)ds\\
&\qquad +p\int_t^T\left(\left(Y^{m}_i(s)-Y^{m}_j(s)-k_{i,j}\right)^{+}\right)^{p-1}\beta^{m}_{i,j}(s)\left(Z^{m}_i(s)-Z^{m}_j(s)\right)ds\\
&\qquad -p\int_t^T\left(\left(Y^{m}_i(s)-Y^{m}_j(s)-k_{i,j}\right)^{+}\right)^{p-1}\left(Z^{m}_i(s)-Z^{m}_j(s)\right)dW(s)\\
&\leq p\gamma\int_t^T\left(\left(Y^{m}_i(s)-Y^{m}_j(s)-k_{i,j}\right)^{+}\right)^{p-1}\left(2+|Y^{m}_i(s)|+|Y^{m}_j(s)|\right)ds\\
&\qquad -p\int_t^T\left(\left(Y^{m}_i(s)-Y^{m}_j(s)-k_{i,j}\right)^{+}\right)^{p-1}\left(Z^{m}_i(s)-Z^{m}_j(s)\right)dW^m_{i,j}(s)
\end{align*}
where $W^m_{i,j}(s)=W(s)-\int_0^s\beta^m_{i,j}(r)dr$ is a Brownian motion under an equivalent probability measure $\frac{dP^m_{i,j}}{dP}=\mathcal{E}_T\left(\beta^m_{i,j}\cdot W\right)$, and the process $\beta^{m}_{i,j}$ satisfies $|\beta^m_{i,j}|\leq\gamma(1+|Z^m_i|+|Z^m_j|)$. Taking conditional expectation with respect to $\mathcal{F}_t$ and $E^{m}_{i,j}$, we have
\begin{align*}
&\left(\left(Y^{m}_i(t)-Y^{m}_j(t)-k_{i,j}\right)^{+}\right)^p+pmE^m_{i,j}\left[\int_t^T\left(\left(Y^{m}_i(s)-Y^{m}_j(s)-k_{i,j}\right)^{+}\right)^pds\bigg|\mathcal{F}_t\right]\\
&\qquad +\frac{p(p-1)}{2}E^m_{i,j}\left[\int_t^T\left(\left(Y^{m}_i(s)-Y^{m}_j(s)-k_{i,j}\right)^{+}\right)^{p-2}|Z^{m}_i(s)-Z^{m}_j(s)|^2ds\bigg|\mathcal{F}_t\right]\\
&\leq p\gamma E^m_{i,j}\left[\int_t^T\left(\left(Y^{m}_i(s)-Y^{m}_j(s)-k_{i,j}\right)^{+}\right)^{p-1}\left(2+|Y^{m}_i(s)|+|Y^{m}_j(s)|\right)ds\bigg|\mathcal{F}_t\right]\\
&\leq
(p-1)mE^m_{i,j}\left[\int_t^T\left(\left(Y^{m}_i(s)-Y^{m}_j(s)-k_{i,j}\right)^{+}\right)^{p}ds\bigg|\mathcal{F}_t\right]\\
&\qquad+\frac{\gamma^p}{m^{p-1}}E^m_{i,j}\left[\int_t^T\left(2+|Y^{m}_i(s)|+|Y^m_j(s)|\right)^pds\bigg|\mathcal{F}_t\right]
\end{align*}
where we have used the Young's inequality in the last inequality. 
Thus, it holds that
\begin{align*}
    \left(\left(Y^{m}_i(t)-Y^{m}_j(t)-k_{i,j}\right)^{+}\right)^p
    \leq\frac{\gamma^p}{m^{p-1}}E^m_{i,j}\left[\int_t^T\left(2+|Y^{m}_i(s)|+|Y^m_j(s)|\right)^pds\bigg|\mathcal{F}_t\right].
    \end{align*}
Noting the uniform boundedness of $Y^m$, we get
\begin{align*}
    \left\|\left(Y^{m}_i-Y^{m}_j-k_{i,j}\right)^{+}\right\|
    \leq\frac{2\gamma\left(1+C_1\right)T^{\frac{1}{p}}}{m^{\frac{p-1}{p}}}.
    \end{align*}
Thefreore, by letting $p\rightarrow +\infty$, there exists a constant $C_3>0$ which is independent of $m$ such that
\begin{equation*}
    \left\|\left(Y^{m}_i(s)-Y^{m}_j(s)-k_{i,j}\right)^{+}\right\|_{\mathcal{S}^{\infty}}\leq \frac{C_3}{m}.
\end{equation*}
\textbf{Step 3.} Recalling that for each $m$ and $t\in[0,T]$, 
\begin{equation*}
    Y^0_i(t)\geq Y^m_i(t)\geq Y^{m+1}(t)\geq\tilde{Y}(t),~~i=1,2,\ldots,n,
\end{equation*}
hence, for a.e. $t$ and P-a.s. $\omega$, $\{Y^m(\omega,t)\}_{m}$ admits a limit which is denoted by $Y(\omega,t)$. Applying Lebesgue's dominated convergence theorem, we have for any $p>2$
\begin{equation*}
    \lim_{m\rightarrow \infty}E\int_0^T|Y^m(s)-Y(s)|^pds=0.
\end{equation*}
In the following we will show that $\{(Y^m,Z^m)\}_{m}$ is a Cauchy sequence in $\mathcal{S}^p_c(\mathbb{R}^n)\times\mathcal{H}^p(\mathbb{R}^{n\times d})$ for any $p>2$. Since it holds that for all $m,\tilde{m}$, $i=1,2,\ldots,n$ and $t\in[0,T]$,
\begin{align*}
    Y^m_i(t)-Y^{\tilde{m}}_i(t)&=\int_t^T\left(h_i(s,Y^m_i(s))-h_i(s,Y^{\tilde{m}}_i(s))+f(s,Z^m_i(s))-f(s,Z^{\tilde{m}}_i(s))\right)ds\\
    &\quad +\int_t^T\left(\tilde{m}\sum_{l=1}^n\left(Y^{\tilde{m}}_i(s)-Y^{\tilde{m}}_l(s)-k_{i,l}\right)^+-m\sum_{l=1}^n\left(Y^{m}_i(s)-Y^{m}_l(s)-k_{i,l}\right)^+\right)ds\\
    &\qquad -\int_t^T\left(Z^{m}_i(s)-Z^{\tilde{m}}_i(s)\right)dW(s)\\
    &=\int_t^T\left(h_i(s,Y^m_i(s))-h_i(s,Y^{\tilde{m}}_i(s))+\beta^{m,\tilde{m}}_i(s)(Z^m_i(s)-Z^{\tilde{m}}_i(s))\right)ds\\
    &\quad +\int_t^T\left(\tilde{m}\sum_{l=1}^n\left(Y^{\tilde{m}}_i(s)-Y^{\tilde{m}}_l(s)-k_{i,l}\right)^+-m\sum_{l=1}^n\left(Y^{m}_i(s)-Y^{m}_l(s)-k_{i,l}\right)^+\right)ds\\
    &\qquad -\int_t^T\left(Z^{m}_i(s)-Z^{\tilde{m}}_i(s)\right)dW(s)\\
    &=\int_t^T\left(h_i(s,Y^m_i(s))-h_i(s,Y^{\tilde{m}}_i(s))\right)ds-\int_t^T\left(Z^{m}_i(s)-Z^{\tilde{m}}_i(s)\right)dW^{m,\tilde{m}}_i(s)\\
    &\quad +\int_t^T\left(\tilde{m}\sum_{l=1}^n\left(Y^{\tilde{m}}_i(s)-Y^{\tilde{m}}_l(s)-k_{i,l}\right)^+-m\sum_{l=1}^n\left(Y^{m}_i(s)-Y^{m}_l(s)-k_{i,l}\right)^+\right)ds
\end{align*}
where $W^{m,\tilde{m}}_{i}(s)=W(s)-\int_0^s\beta^{m,\tilde{m}}_{i}(r)dr$ is a Brownian motion under an equivalent probability measure $\frac{dP^{m,\tilde{m}}_{i}}{dP}=\mathcal{E}_T\left(\beta^{m,\tilde{m}}_{i}\cdot W\right)$, and the process $\beta^{m,\tilde{m}}_{i}$ satisfies $|\beta^{m,\tilde{m}}_{i}|\leq\gamma(1+|Z^m_i|+|Z^{\tilde{m}}_i|)$.

For any $p>2$, applying It\^{o}'s formula to $|Y^m_i(t)-Y^{\tilde{m}}_i(t)|^p$ and taking conditional expectation with respect to $\mathcal{F}_t$ and $E^{m,\tilde{m}}_i$, noting the uniform boundedness of $Y^m$, we get
\begin{align*}
    &|Y^m_i(t)-Y^{\tilde{m}}_i(t)|^p+\frac{p(p-1)}{2}E^{m,\tilde{m}}_i\left[\int_t^T|Y^m_i(s)-Y^{\tilde{m}}_i(s)|^{p-2}\left(Z^{m}_i(s)-Z^{\tilde{m}}_i(s)\right)^2ds\bigg|\mathcal{F}_t\right]\\
    &=pE^{m,\tilde{m}}_i\left[\int_t^T|Y^m_i(s)-Y^{\tilde{m}}_i(s)|^{p-1}\left(h_i(s,Y^m_i(s))-h_i(s,Y^{\tilde{m}}_i(s))\right)ds\bigg|\mathcal{F}_t\right]\\
    &+pE^{m,\tilde{m}}_i\left[\int_t^T|Y^m_i(s)-Y^{\tilde{m}}_i(s)|^{p-1}\left(\tilde{m}\sum_{l=1}^n\left(Y^{\tilde{m}}_i(s)-Y^{\tilde{m}}_l(s)-k_{i,l}\right)^+-m\sum_{l=1}^n\left(Y^{m}_i(s)-Y^{m}_l(s)-k_{i,l}\right)^+\right)ds\bigg|\mathcal{F}_t\right]\\
    &\leq pE^{m,\tilde{m}}_i\left[\int_t^T|Y^m_i(s)-Y^{\tilde{m}}_i(s)|^{p-1}\left|h_i(s,Y^m_i(s))-h_i(s,Y^{\tilde{m}}_i(s))\right|ds\bigg|\mathcal{F}_t\right]\\
    &+p\sum_{l=1}^nE^{m,\tilde{m}}_i\left[\int_t^T|Y^m_i(s)-Y^{\tilde{m}}_i(s)|^{p-1}\left|\tilde{m}\left(Y^{\tilde{m}}_i(s)-Y^{\tilde{m}}_l(s)-k_{i,l}\right)^+-m\left(Y^{m}_i(s)-Y^{m}_l(s)-k_{i,l}\right)^+\right|ds\bigg|\mathcal{F}_t\right]\\
    &\leq 2\gamma p(1+C_1)E^{m,\tilde{m}}_i\left[\int_t^T|Y^m_i(s)-Y^{\tilde{m}}_i(s)|^{p-1}ds\bigg|\mathcal{F}_t\right]\\
    &+p\sum_{l=1}^nE^{m,\tilde{m}}_i\left[\int_t^T|Y^m_i(s)-Y^{\tilde{m}}_i(s)|^{p-1}\left|\tilde{m}\left(Y^{\tilde{m}}_i(s)-Y^{\tilde{m}}_l(s)-k_{i,l}\right)^+\right|ds\bigg|\mathcal{F}_t\right]\\
    &+p\sum_{l=1}^nE^{m,\tilde{m}}_i\left[\int_t^T|Y^m_i(s)-Y^{\tilde{m}}_i(s)|^{p-1}\left|m\left(Y^{m}_i(s)-Y^{m}_l(s)-k_{i,l}\right)^+\right|ds\bigg|\mathcal{F}_t\right].
\end{align*}
Therefore, it holds that
\begin{align*}
    &\sup_{0\leq t\leq T}|Y^m_i(t)-Y^{\tilde{m}}_i(t)|^p \\
    &\leq 2\gamma p(1+C_1)\sup_{0\leq t\leq T}E^{m,\tilde{m}}_i\left[\int_0^T|Y^m_i(s)-Y^{\tilde{m}}_i(s)|^{p-1}ds\bigg|\mathcal{F}_t\right]\\
    &+p\sum_{l=1}^n\sup_{0\leq t\leq T}E^{m,\tilde{m}}_i\left[\int_0^T|Y^m_i(s)-Y^{\tilde{m}}_i(s)|^{p-1}\left|\tilde{m}\left(Y^{\tilde{m}}_i(s)-Y^{\tilde{m}}_l(s)-k_{i,l}\right)^+\right|ds\bigg|\mathcal{F}_t\right]\\
    &+p\sum_{l=1}^n\sup_{0\leq t\leq T}E^{m,\tilde{m}}_i\left[\int_0^T|Y^m_i(s)-Y^{\tilde{m}}_i(s)|^{p-1}\left|m\left(Y^{m}_i(s)-Y^{m}_l(s)-k_{i,l}\right)^+\right|ds\bigg|\mathcal{F}_t\right].
\end{align*}
Hence, it follows from Doob's maximal inequality that
\begin{align*}
    &E^{m,\tilde{m}}_i\left[\sup_{0\leq t\leq T}|Y^m_i(t)-Y^{\tilde{m}}_i(t)|^p\right]\\
    &\leq 4\gamma p(1+C_1)E^{m,\tilde{m}}_i\left[\left(\int_0^T|Y^m_i(s)-Y^{\tilde{m}}_i(s)|^{p-1}ds\right)^2\right]^{\frac{1}{2}}\\
    &+2p\sum_{l=1}^nE^{m,\tilde{m}}_i\left[\left(\int_0^T|Y^m_i(s)-Y^{\tilde{m}}_i(s)|^{p-1}\left|\tilde{m}\left(Y^{\tilde{m}}_i(s)-Y^{\tilde{m}}_l(s)-k_{i,l}\right)^+\right|ds\right)^2\right]^{\frac{1}{2}}\\
    &+2p\sum_{l=1}^nE^{m,\tilde{m}}_i\left[\left(\int_0^T|Y^m_i(s)-Y^{\tilde{m}}_i(s)|^{p-1}\left|m\left(Y^{m}_i(s)-Y^{m}_l(s)-k_{i,l}\right)^+\right|ds\right)^2\right]^{\frac{1}{2}}.
\end{align*}
Since $\|\beta^{m,\tilde{m}}_{i}\cdot W\|_{BMO}$ is independent of $m$ and $\frac{dP^{m,\tilde{m}}_{i}}{dP}=\mathcal{E}_T\left(\beta^{m,\tilde{m}}_{i}\cdot W\right)$, there exists $p_2>1$ and $p_3>1$ such that
\begin{equation*}
    E^{m,\tilde{m}}_{i}\left[\left(\frac{dP}{dP^{m,\tilde{m}}_{i}}\right)^{p_2}\right]<+\infty,~~~E\left[\left(\frac{dP^{m,\tilde{m}}_i}{dP}\right)^{p_3}\right]<+\infty.
\end{equation*}
Thus, we obtain
\begin{align*}
    &E^{m,\tilde{m}}_i\left[\sup_{0\leq t\leq T}|Y^m_i(t)-Y^{\tilde{m}}_i(t)|^p\right]\\
    &\leq 4\gamma p(1+C_1)E\left[\left(\frac{dP^{m,\tilde{m}}_i}{dP}\right)^{p_2}\right]^{\frac{1}{2p_2}}E\left[\left(\int_0^T|Y^m_i(s)-Y^{\tilde{m}}_i(s)|^{p-1}ds\right)^{2q_2}\right]^{\frac{1}{2q_2}}\\
    &+2p\sum_{l=1}^nE\left[\left(\frac{dP^{m,\tilde{m}}_i}{dP}\right)^{p_2}\right]^{\frac{1}{2p_2}}E\left[\left(\int_0^T|Y^m_i(s)-Y^{\tilde{m}}_i(s)|^{p-1}\left|\tilde{m}\left(Y^{\tilde{m}}_i(s)-Y^{\tilde{m}}_l(s)-k_{i,l}\right)^+\right|ds\right)^{2q_2}\right]^{\frac{1}{2q_2}}\\
    &+2p\sum_{l=1}^nE\left[\left(\frac{dP^{m,\tilde{m}}_i}{dP}\right)^{p_2}\right]^{\frac{1}{2p_2}}E\left[\left(\int_0^T|Y^m_i(s)-Y^{\tilde{m}}_i(s)|^{p-1}\left|m\left(Y^{m}_i(s)-Y^{m}_l(s)-k_{i,l}\right)^+\right|ds\right)^{2q_2}\right]^{\frac{1}{2q_2}}\\
    &\leq 4\gamma p(1+C_1)TE\left[\left(\frac{dP^{m,\tilde{m}}_i}{dP}\right)^{p_2}\right]^{\frac{1}{2p_2}}E\left[\int_0^T|Y^m_i(s)-Y^{\tilde{m}}_i(s)|^{2q_2(p-1)}ds\right]^{\frac{1}{2q_2}}\\
    &+2p\sum_{l=1}^nE\left[\left(\frac{dP^{m,\tilde{m}}_i}{dP}\right)^{p_2}\right]^{\frac{1}{2p_2}}E\left[\int_0^T|Y^m_i(s)-Y^{\tilde{m}}_i(s)|^{4q_2(p-1)}ds\right]^{\frac{1}{4q_2}}E\left[\int_0^T\left|\tilde{m}\left(Y^{\tilde{m}}_i(s)-Y^{\tilde{m}}_l(s)-k_{i,l}\right)^+\right|^{4q_2}ds\right]^{\frac{1}{4q_2}}\\
    &+2p\sum_{l=1}^nE\left[\left(\frac{dP^{m,\tilde{m}}_i}{dP}\right)^{p_2}\right]^{\frac{1}{2p_2}}E\left[\int_0^T|Y^m_i(s)-Y^{\tilde{m}}_i(s)|^{4q_2(p-1)}ds\right]^{\frac{1}{4q_2}}E\left[\int_0^T\left|m\left(Y^{m}_i(s)-Y^{m}_l(s)-k_{i,l}\right)^+\right|^{4q_2}ds\right]^{\frac{1}{4q_2}}
\end{align*}
where we have used H\"{o}lder's inequality, Jensen's inequality and $\frac{1}{p_2}+\frac{1}{q_2}=1$. On the other hand, it holds that
\begin{align*}
    &E\left[\sup_{0\leq t\leq T}|Y^m_i(t)-Y^{\tilde{m}}_i(t)|^p\right]\\
    &\leq E^{m,\tilde{m}}_i\left[\frac{dP}{dP^{m,\tilde{m}}_i}\sup_{0\leq t\leq T}|Y^m_i(t)-Y^{\tilde{m}}_i(t)|^p\right]\\
    &\leq E^{m,\tilde{m}}_i\left[\left(\frac{dP}{dP^{m,\tilde{m}}_i}\right)^{p_3}\right]^{\frac{1}{p_3}}E^{m,\tilde{m}}_i\left[\sup_{0\leq t\leq T}|Y^m_i(t)-Y^{\tilde{m}}_i(t)|^{pq_3}\right]^{\frac{1}{q_3}}
\end{align*}
where $\frac{1}{p_3}+\frac{1}{q_3}=1$.

Concluding above, there exists a constant $C_4>0$ (independent of m) such that for any $p>2$
\begin{align*}
    &E\left[\sup_{0\leq t\leq T}|Y^m_i(t)-Y^{\tilde{m}}_i(t)|^p\right]\\
    &\leq C_4 E\left[\int_0^T|Y^m_i(s)-Y^{\tilde{m}}_i(s)|^{2q_2(pq_3-1)}ds\right]^{\frac{1}{2q_2q_3}} + C_4 E\left[\int_0^T|Y^m_i(s)-Y^{\tilde{m}}_i(s)|^{4q_2(pq_3-1)}ds\right]^{\frac{1}{4q_2q_3}}
\end{align*}
which implies that
\begin{equation*}
    \lim_{m,\tilde{m}\rightarrow +\infty} E\left[\sup_{0\leq t\leq T}|Y^m_i(t)-Y^{\tilde{m}}_i(t)|^p\right]=0.
\end{equation*}
Therefore, we deduce that $\{Y^m\}_m$ is a Cauchy sequence in $\mathcal{S}^p_c(\mathbb{R}^n)$ for any $p>2$. In particular, it holds that
\begin{equation*}
    \lim_{m\rightarrow +\infty}\|Y^m-Y\|_{\mathcal{S}^{p}}=0.
\end{equation*}
We will now show that $\{Z^m\}_m$ is a Cauchy sequence in $\mathcal{H}^p(\mathbb{R}^{n\times d})$ for any $p>2$. In fact, applying It\^{o}'s formula to $\left(Y^m_i(t)-Y^{\tilde{m}}_i(t)\right)^2$, we get
\begin{align*}
    &\left(Y^m_i(t)-Y^{\tilde{m}}_i(t)\right)^2+\int_t^T\left(Z^{m}_i(s)-Z^{\tilde{m}}_i(s)\right)^2ds\\
    &=2\int_t^T\left(Y^m_i(s)-Y^{\tilde{m}}_i(s)\right)\left(h_i(s,Y^m_i(s))-h_i(s,Y^{\tilde{m}}_i(s))\right)ds\\
    &+2\int_t^T\left(Y^m_i(s)-Y^{\tilde{m}}_i(s)\right)\left(\tilde{m}\sum_{l=1}^n\left(Y^{\tilde{m}}_i(s)-Y^{\tilde{m}}_l(s)-k_{i,l}\right)^+-m\sum_{l=1}^n\left(Y^{m}_i(s)-Y^{m}_l(s)-k_{i,l}\right)^+\right)ds\\
    &-2\int_t^T\left(Y^m_i(s)-Y^{\tilde{m}}_i(s)\right)\left(Z^{m}_i(s)-Z^{\tilde{m}}_i(s)\right)dW^{m,\tilde{m}}_i(s)
\end{align*}
which implies that
\begin{align*}
    &\int_t^T\left(Z^{m}_i(s)-Z^{\tilde{m}}_i(s)\right)^2ds\\
    &\leq 2\int_t^T\left|Y^m_i(s)-Y^{\tilde{m}}_i(s)\right|\left|h_i(s,Y^m_i(s))-h_i(s,Y^{\tilde{m}}_i(s))\right|ds\\
    &+2\sum_{l=1}^n\int_t^T\left|Y^m_i(s)-Y^{\tilde{m}}_i(s)\right|\left|\tilde{m}\left(Y^{\tilde{m}}_i(s)-Y^{\tilde{m}}_l(s)-k_{i,l}\right)^+-m\left(Y^{m}_i(s)-Y^{m}_l(s)-k_{i,l}\right)^+\right|ds\\
    &+2\left|\int_t^T\left(Y^m_i(s)-Y^{\tilde{m}}_i(s)\right)\left(Z^{m}_i(s)-Z^{\tilde{m}}_i(s)\right)dW^{m,\tilde{m}}_i(s)\right|.
\end{align*}
Thus, it holds that for any $p>2$
\begin{align*}
    &\left(\int_t^T\left(Z^{m}_i(s)-Z^{\tilde{m}}_i(s)\right)^2ds\right)^{\frac{p}{2}}\\
    &\leq 2^{\frac{p}{2}}(n+2)^{\frac{p}{2}}\left(\int_t^T\left|Y^m_i(s)-Y^{\tilde{m}}_i(s)\right|\left|h_i(s,Y^m_i(s))-h_i(s,Y^{\tilde{m}}_i(s))\right|ds\right)^{\frac{p}{2}}\\
    &+2^{\frac{p}{2}}(n+2)^{\frac{p}{2}}\sum_{l=1}^n\left(\int_t^T\left|Y^m_i(s)-Y^{\tilde{m}}_i(s)\right|\left|\tilde{m}\left(Y^{\tilde{m}}_i(s)-Y^{\tilde{m}}_l(s)-k_{i,l}\right)^+-m\left(Y^{m}_i(s)-Y^{m}_l(s)-k_{i,l}\right)^+\right|ds\right)^{\frac{p}{2}}\\
    &+2^{\frac{p}{2}}(n+2)^{\frac{p}{2}}\left|\int_t^T\left(Y^m_i(s)-Y^{\tilde{m}}_i(s)\right)\left(Z^{m}_i(s)-Z^{\tilde{m}}_i(s)\right)dW^{m,\tilde{m}}_i(s)\right|^{\frac{p}{2}}.
\end{align*}
Therefore, using Burkholder-Davis-Gundy inequality, there exists a constant $C_p>0$ only depending on $p$ such that
\begin{align*}
    &E^{m,\tilde{m}}_i\left[\left(\int_0^T\left(Z^{m}_i(s)-Z^{\tilde{m}}_i(s)\right)^2ds\right)^{\frac{p}{2}}\right]\\
    &\leq 2^{\frac{p}{2}}(n+2)^{\frac{p}{2}}E^{m,\tilde{m}}_i\left[\left(\int_0^T\left|Y^m_i(s)-Y^{\tilde{m}}_i(s)\right|\left|h_i(s,Y^m_i(s))-h_i(s,Y^{\tilde{m}}_i(s))\right|ds\right)^{\frac{p}{2}}\right]\\
    &+2^{\frac{p}{2}}(n+2)^{\frac{p}{2}}\sum_{l=1}^nE^{m,\tilde{m}}_i\left[\left(\int_0^T\left|Y^m_i(s)-Y^{\tilde{m}}_i(s)\right|\left|\tilde{m}\left(Y^{\tilde{m}}_i(s)-Y^{\tilde{m}}_l(s)-k_{i,l}\right)^+-m\left(Y^{m}_i(s)-Y^{m}_l(s)-k_{i,l}\right)^+\right|ds\right)^{\frac{p}{2}}\right]\\
    &+2^{\frac{p}{2}}(n+2)^{\frac{p}{2}}E^{m,\tilde{m}}_i\left[\left|\int_0^T\left(Y^m_i(s)-Y^{\tilde{m}}_i(s)\right)\left(Z^{m}_i(s)-Z^{\tilde{m}}_i(s)\right)dW^{m,\tilde{m}}_i(s)\right|^{\frac{p}{2}}\right]\\
    &\leq 2^{\frac{p}{2}}(n+2)^{\frac{p}{2}}E^{m,\tilde{m}}_i\left[\left(\int_0^T\left|Y^m_i(s)-Y^{\tilde{m}}_i(s)\right|\left|h_i(s,Y^m_i(s))-h_i(s,Y^{\tilde{m}}_i(s))\right|ds\right)^{\frac{p}{2}}\right]\\
    &+2^{\frac{p}{2}}(n+2)^{\frac{p}{2}}\sum_{l=1}^nE^{m,\tilde{m}}_i\left[\left(\int_0^T\left|Y^m_i(s)-Y^{\tilde{m}}_i(s)\right|\left|\tilde{m}\left(Y^{\tilde{m}}_i(s)-Y^{\tilde{m}}_l(s)-k_{i,l}\right)^+-m\left(Y^{m}_i(s)-Y^{m}_l(s)-k_{i,l}\right)^+\right|ds\right)^{\frac{p}{2}}\right]\\
    &+2^{\frac{p}{2}}(n+2)^{\frac{p}{2}}C_{p}E^{m,\tilde{m}}_i\left[\left(\int_0^T\left(Y^m_i(s)-Y^{\tilde{m}}_i(s)\right)^2\left(Z^{m}_i(s)-Z^{\tilde{m}}_i(s)\right)^2d(s)\right)^{\frac{p}{4}}\right]\\
    &\leq 2^{\frac{p}{2}}(n+2)^{\frac{p}{2}}E^{m,\tilde{m}}_i\left[\left(\sup_{0\leq s\leq T}\left|Y^m_i(s)-Y^{\tilde{m}}_i(s)\right|\int_0^T\left|h_i(s,Y^m_i(s))-h_i(s,Y^{\tilde{m}}_i(s))\right|ds\right)^{\frac{p}{2}}\right]\\
    &+2^{p}(n+2)^{\frac{p}{2}}\sum_{l=1}^nE^{m,\tilde{m}}_i\left[\left(\sup_{0\leq s\leq T}\left|Y^m_i(s)-Y^{\tilde{m}}_i(s)\right|\int_0^T\left|\tilde{m}\left(Y^{\tilde{m}}_i(s)-Y^{\tilde{m}}_l(s)-k_{i,l}\right)^+\right|ds\right)^{\frac{p}{2}}\right]\\
    &+2^{p}(n+2)^{\frac{p}{2}}\sum_{l=1}^nE^{m,\tilde{m}}_i\left[\left(\sup_{0\leq s\leq T}\left|Y^m_i(s)-Y^{\tilde{m}}_i(s)\right|\int_0^T\left|m\left(Y^{m}_i(s)-Y^{m}_l(s)-k_{i,l}\right)^+\right|ds\right)^{\frac{p}{2}}\right]\\
    &+2^{\frac{p}{2}}(n+2)^{\frac{p}{2}}C_{p}E^{m,\tilde{m}}_i\left[\left(\sup_{0\leq s\leq T}\left(Y^m_i(s)-Y^{\tilde{m}}_i(s)\right)^2\int_0^T\left(Z^{m}_i(s)-Z^{\tilde{m}}_i(s)\right)^2d(s)\right)^{\frac{p}{4}}\right]\\
    &\leq 2^{\frac{p}{2}}(n+2)^{\frac{p}{2}}E^{m,\tilde{m}}_i\left[\left(\sup_{0\leq s\leq T}\left|Y^m_i(s)-Y^{\tilde{m}}_i(s)\right|\int_0^T\left|h_i(s,Y^m_i(s))-h_i(s,Y^{\tilde{m}}_i(s))\right|ds\right)^{\frac{p}{2}}\right]\\
    &+2^{p}(n+2)^{\frac{p}{2}}\sum_{l=1}^nE^{m,\tilde{m}}_i\left[\left(\sup_{0\leq s\leq T}\left|Y^m_i(s)-Y^{\tilde{m}}_i(s)\right|\int_0^T\left|\tilde{m}\left(Y^{\tilde{m}}_i(s)-Y^{\tilde{m}}_l(s)-k_{i,l}\right)^+\right|ds\right)^{\frac{p}{2}}\right]\\
    &+2^{p}(n+2)^{\frac{p}{2}}\sum_{l=1}^nE^{m,\tilde{m}}_i\left[\left(\sup_{0\leq s\leq T}\left|Y^m_i(s)-Y^{\tilde{m}}_i(s)\right|\int_0^T\left|m\left(Y^{m}_i(s)-Y^{m}_l(s)-k_{i,l}\right)^+\right|ds\right)^{\frac{p}{2}}\right]\\
    &+2^{p-1}(n+2)^{p}C^2_{p}E^{m,\tilde{m}}_i\left[\sup_{0\leq s\leq T}\left(Y^m_i(s)-Y^{\tilde{m}}_i(s)\right)^p\right]+\frac{1}{2}E^{m,\tilde{m}}_i\left[\left(\int_0^T\left(Z^{m}_i(s)-Z^{\tilde{m}}_i(s)\right)^2d(s)\right)^{\frac{p}{2}}\right]
\end{align*}
where we have used Young's inequality in the last inequality. Hence, we obtain
\begin{align*}
    &\frac{1}{2}E^{m,\tilde{m}}_i\left[\left(\int_0^T\left(Z^{m}_i(s)-Z^{\tilde{m}}_i(s)\right)^2ds\right)^{\frac{p}{2}}\right]\\
    &\leq 2^{\frac{p}{2}}(n+2)^{\frac{p}{2}}E^{m,\tilde{m}}_i\left[\left(\sup_{0\leq s\leq T}\left|Y^m_i(s)-Y^{\tilde{m}}_i(s)\right|\int_0^T\left|h_i(s,Y^m_i(s))-h_i(s,Y^{\tilde{m}}_i(s))\right|ds\right)^{\frac{p}{2}}\right]\\
    &+2^{p}(n+2)^{\frac{p}{2}}\sum_{l=1}^nE^{m,\tilde{m}}_i\left[\left(\sup_{0\leq s\leq T}\left|Y^m_i(s)-Y^{\tilde{m}}_i(s)\right|\int_0^T\left|\tilde{m}\left(Y^{\tilde{m}}_i(s)-Y^{\tilde{m}}_l(s)-k_{i,l}\right)^+\right|ds\right)^{\frac{p}{2}}\right]\\
    &+2^{p}(n+2)^{\frac{p}{2}}\sum_{l=1}^nE^{m,\tilde{m}}_i\left[\left(\sup_{0\leq s\leq T}\left|Y^m_i(s)-Y^{\tilde{m}}_i(s)\right|\int_0^T\left|m\left(Y^{m}_i(s)-Y^{m}_l(s)-k_{i,l}\right)^+\right|ds\right)^{\frac{p}{2}}\right]\\
    &+2^{p-1}(n+2)^{p}C^2_{p}E^{m,\tilde{m}}_i\left[\sup_{0\leq s\leq T}\left(Y^m_i(s)-Y^{\tilde{m}}_i(s)\right)^p\right]\\
    &\leq 2^{\frac{p}{2}}(n+2)^{\frac{p}{2}}E\left[\left(\frac{dP^{m,\tilde{m}_i}}{dP}\right)^{p_2}\right]^{\frac{1}{p_2}}E\left[\left(\sup_{0\leq s\leq T}\left|Y^m_i(s)-Y^{\tilde{m}}_i(s)\right|\int_0^T\left|h_i(s,Y^m_i(s))-h_i(s,Y^{\tilde{m}}_i(s))\right|ds\right)^{\frac{pq_2}{2}}\right]^{\frac{1}{q_2}}\\
    &+2^{p}(n+2)^{\frac{p}{2}}\sum_{l=1}^nE\left[\left(\frac{dP^{m,\tilde{m}_i}}{dP}\right)^{p_2}\right]^{\frac{1}{p_2}}\left[\left(\sup_{0\leq s\leq T}\left|Y^m_i(s)-Y^{\tilde{m}}_i(s)\right|\int_0^T\left|\tilde{m}\left(Y^{\tilde{m}}_i(s)-Y^{\tilde{m}}_l(s)-k_{i,l}\right)^+\right|ds\right)^{\frac{pq_2}{2}}\right]^{\frac{1}{q_2}}\\
    &+2^{p}(n+2)^{\frac{p}{2}}\sum_{l=1}^nE\left[\left(\frac{dP^{m,\tilde{m}_i}}{dP}\right)^{p_2}\right]^{\frac{1}{p_2}}E\left[\left(\sup_{0\leq s\leq T}\left|Y^m_i(s)-Y^{\tilde{m}}_i(s)\right|\int_0^T\left|m\left(Y^{m}_i(s)-Y^{m}_l(s)-k_{i,l}\right)^+\right|ds\right)^{\frac{pq_2}{2}}\right]^{\frac{1}{q_2}}\\
    &+2^{p-1}(n+2)^{p}C^2_{p}E\left[\left(\frac{dP^{m,\tilde{m}_i}}{dP}\right)^{p_2}\right]^{\frac{1}{p_2}}E\left[\sup_{0\leq s\leq T}\left(Y^m_i(s)-Y^{\tilde{m}}_i(s)\right)^{pq_2}\right]^{\frac{1}{q_2}}.
\end{align*}
On the other hand, it holds that for any $p>2$
\begin{align*}
    E\left[\left(\int_0^T\left(Z^{m}_i(s)-Z^{\tilde{m}}_i(s)\right)^2ds\right)^{\frac{p}{2}}\right]&=E^{m,\tilde{m}}_i\left[\frac{dP}{dP^{m,\tilde{m}}_i}\left(\int_0^T\left(Z^{m}_i(s)-Z^{\tilde{m}}_i(s)\right)^2ds\right)^{\frac{p}{2}}\right]\\
    &\leq E^{m,\tilde{m}}_i\left[\left(\frac{dP}{dP^{m,\tilde{m}}_i}\right)^{p_3}\right]^{\frac{1}{p_3}}E^{m,\tilde{m}}_i\left[\left(\int_0^T\left(Z^{m}_i(s)-Z^{\tilde{m}}_i(s)\right)^2ds\right)^{\frac{pq_3}{2}}\right]^{\frac{1}{q_3}}.
\end{align*}
Combining the above, for any $p>2$, we get
\begin{equation*}
    \lim_{m,\tilde{m}\rightarrow +\infty}\|Z^m-Z^{\tilde{m}}\|_{\mathcal{H}^p}=0
\end{equation*}
which implies $\{Z^m\}_{m}$ is a Cauchy sequence in $\mathcal{H}^p(\mathbb{R}^{n\times d})$ for any $p>2$. In particular, there exists $Z\in \mathcal{H}^{p}(\mathbb{R}^{n\times d})$ for any $p>2$ such that
\begin{equation*}
    \lim_{m\rightarrow +\infty}\|Z^m-Z\|_{\mathcal{H}^p}=0.
\end{equation*}
Now we define the non-decreasing process $K^m_i$ as follows
\begin{equation*}
    K^m_i(t)=m\int_0^t\sum_{l=1}^n\left(Y^m_i(s)-Y^m_l(s)-k_{i,l}\right)^+ds,~~t\in[0,T],i=1,2,\ldots,n.
\end{equation*}
Thus it follows from the penalized BSDE \eqref{BSDE_penalized} that
\begin{equation*}
    K^m_i(t)=Y^m_i(t)-Y^m_i(0)+\int_0^tg_i(s,Y^m_i(s),Z^m_i(s))ds-\int_0^tZ^m_i(s)dW(s),~~~i=1,2,\ldots,n.
\end{equation*}
Setting
\begin{equation*}
    K_i(t)=Y_i(t)-Y_i(0)+\int_0^tg_i(s,Y_i(s),Z_i(s))ds-\int_0^tZ_i(s)dW(s),~~~i=1,2,\ldots,n,
\end{equation*}
and noting that $\lim_{m\rightarrow +\infty}\|Y^m-Y\|_{\mathcal{S}^{p}}=0$ and $\lim_{m\rightarrow +\infty}\|Z^m-Z\|_{\mathcal{H}^p}=0$ for any $p>2$, one can easily check that
\begin{equation*}
    \lim_{m\rightarrow +\infty}\|K^m-K\|_{\mathcal{S}^p}=0,
\end{equation*}
which implies that $\left(Y,Z,K\right)$ satisfies the first relation in RBSDE \eqref{BSDE}.

\noindent \textbf{Step 4.} We will now prove that $\left(Y,Z,K\right)$ satisfies second and third realtions in RBSDE \eqref{BSDE}. In fact, noting that $\lim_{m\rightarrow +\infty}\|Y^m-Y\|_{\mathcal{S}^{p}}=0$ for any $p>2$ and letting $m\rightarrow +\infty$ in \eqref{eq:uniform_penalized}, we deduce that for any $p>2$
\begin{equation*}
    E\left[\sup_{0\leq t\leq T}\left(\left(Y_i(t)-Y_j(t)-k_{i,j}\right)^+\right)^p\right]=0,~~~i,j=1,2,\ldots,n.
\end{equation*}
Thus
\begin{equation*}
    Y_i(t)\leq Y_j(t)+k_{i,j},~~t\in[0,T],~~i,j=1,2,\ldots,n
\end{equation*}
which is equivalent to the fact the for $P$-a.s. $\omega$, 
\begin{equation*}
    Y(\omega,t)\in\overline{Q},~~~\forall t\in[0,T].
\end{equation*}
On the other hand, for $t\in[0,T]$ and $i,j,l=1,2,\ldots,n$, it holds that 
\begin{equation*}
    \left(Y^m_i(t)-Y^m_l(t)-k_{i,l}\right)^-\left(Y^m_i(t)-Y^m_l(t)-k_{i,l}\right)^+=0  
\end{equation*}
and if $l=i$
\begin{align*}
\left(Y^m_i(t)-Y^m_j(t)-k_{i,j}\right)^-\left(Y^m_i(t)-Y^m_l(t)-k_{i,l}\right)^+=0.
\end{align*}
Thus, we get
\begin{equation*} 
    \min_{j\neq i}\left\{\left(Y^m_i(t)-Y^m_j(t)-k_{i,j}\right)^-\left(Y^m_i(t)-Y^m_l(t)-k_{i,l}\right)^+\right\}=0.
\end{equation*}
Therefore, it holds that for $i=1,2,\ldots,n$
\begin{align*}
    &\int_0^T\left(Y^m_i(t)-\min_{j\neq i}\left(Y^m_j(t)+k_{i,j}\right)\right)^-dK^m_i(t)\\
    &=m\sum_{l=1}^n\int_0^T\min_{j\neq i}\left\{\left(Y^m_i(t)-Y^m_j(t)-k_{i,j}\right)^-\left(Y^m_i(t)-Y^m_l(t)-k_{i,l}\right)^+\right\}dt=0.
\end{align*}
Letting $m\rightarrow +\infty$ and applying \cite[Lemma 5.8]{GP}, we obtain
\begin{equation*}
    \int_0^T\left(Y_i(t)-\min_{j\neq i}\left(Y_j(t)+k_{i,j}\right)\right)^-dK_i(t)=0,~~~i=1,2,\ldots,n.  
\end{equation*}
Hence, $(Y,Z,K)$ is a solution of RBSDE \eqref{BSDE} satisfying $(Y,Z,K)\in\mathcal{S}^p_c(\mathbb{R}^n)\times\mathcal{H}^p(\mathbb{R}^{n\times d})\times \mathcal{A}^p_c(\mathbb{R}^n)$ for any $p>2$.

\noindent \textbf{Step 5.} Finally, we show that $(Y,Z\cdot W,K)\in\mathcal{S}^{\infty}_c(\mathbb{R}^n)\times BMO\times \mathcal{A}^{\infty}_c(\mathbb{R}^n)$. Indeed, it follows from the uniform boundedness of $Y^m$ and $\lim_{m\rightarrow +\infty}\|Y^m-Y\|_{\mathcal{S}^{\infty}}=0$ for any $p>2$ that $Y\in\mathcal{S}^{\infty}_c(\mathbb{R}^n)$. Similarly, we get $K\in\mathcal{A}^{\infty}_c(\mathbb{R}^n)$. 

For any $i=1,2,\ldots,n$ and all $m$, applying It\^{o} formula to $e^{2(\gamma+1) Y_i(t)}$ implies that for all $t\in[0,T]$
\begin{align*}
    e^{2(\gamma+1)Y_i(t)}&=e^{2(\gamma+1)\xi_i}-\int_{t}^T2(\gamma+1)e^{2(\gamma+1)Y_i(s)}Z_i(s)dW(s)--\int_{t}^T2(\gamma+1)e^{2(\gamma+1)Y_i(s)}dK_i(s)\\
    &+\int_t^T\left(2(\gamma+1)e^{2(\gamma+1)Y_i(s)}g_i\left(s,Y_i(s),Z_i(s)\right)-2(\gamma+1)^2e^{2(\gamma+1)Y_i(s)}|Z_i(s)|^2\right)ds\\
    &\leq e^{2(\gamma+1)\xi_i}-\int_{t}^T2(\gamma+1)e^{2(\gamma+1)Y_i(s)}Z_i(s)dW(s)\\
    &+\int_t^T\left(2\gamma(\gamma+1)e^{2(\gamma+1)Y_i(s)}\left(2+|Y_i(s)|+|Z_i(s)|^2\right)-2(\gamma+1)^2e^{2(\gamma+1)Y_i(s)}|Z_i(s)|^2\right)ds\\
    &\leq e^{2(\gamma+1)\xi_i}-\int_{t}^T2(\gamma+1)e^{2(\gamma+1)Y_i(s)}Z_i(s)dW(s)\\
    &+\int_t^T\left(2\gamma(\gamma+1)e^{2(\gamma+1)Y_i(s)}\left(2+|Y_i(s)|\right)-2(\gamma+1)e^{2(\gamma+1)Y_i(s)}|Z_i(s)|^2\right)ds.
\end{align*}
Hence, it holds that
\begin{align*}
    2(\gamma+1)\int_t^Te^{2(\gamma+1)Y_i(s)}|Z_i(s)|^2ds&\leq e^{2(\gamma+1)\xi_i}-\int_{t}^T2(\gamma+1)e^{2(\gamma+1)Y_i(s)}Z_i(s)dW(s)\\
    &+\int_t^T2\gamma(\gamma+1)e^{2(\gamma+1)Y_i(s)}\left(2+|Y_i(s)|\right)ds.
\end{align*}
Taking conditional expectation with respect to $\mathcal{F}_t$ and noting the boundedness of $Y$, it holds that
\begin{align*}
    E\left[\int_t^T|Z_i(s)|^2ds\bigg|\mathcal{F}_t\right]\leq \frac{1}{2}e^{2(\gamma+1)\|Y_i\|_{\mathcal{S}^{\infty}}}\left(e^{2(\gamma+1)\|\xi_i\|_{\infty}}+2\gamma(\gamma+1)e^{2(\gamma+1)\|Y_i\|_{\mathcal{S}^{\infty}}}\left(2+\|Y_i\|_{\mathcal{S}^{\infty}}\right)T\right)
\end{align*}
which implies that $Z\cdot W\in BMO$.
\end{proof}
\subsection{Uniqueness}
Following the idea of \cite{HT1}, we prove the uniqueness by verifying that $Y$ is the value of an optimal switching problem of quadratic BSDEs, for a solution $(Y,Z,K)$ of RBSDE \eqref{BSDE} satisfying $(Y,Z\cdot W,K)\in\mathcal{S}^{\infty}_c(\mathbb{R}^n)\times BMO \times\mathcal{A}^{\infty}_c(\mathbb{R}^n)$. We now introduce the optimal switching problem.
\begin{definition}
    Let $\{\tau_j\}^{\infty}_{j=0}$ be a non-decreasing sequence of stopping times with values in $[0,T]$ and for any $j$, $\alpha_j$ is an $\mathcal{F}_{\tau_j}$-measurable random variable taking values in $\{1,2,\ldots,n\}$. By defining 
    \begin{equation*}
        a(t):=a_0\mathbf{1}_{\{\tau_0\}}(t)+\sum_{j=1}^{\infty}\alpha_{j-1}\mathbf{1}_{(\tau_{j-1},\tau_j]}(t),~~t\in[0,T],
    \end{equation*}
     we call the sequence $\{(\tau_j,\alpha_j)\}$ or $a$ an admissible switching strategy starting from mode $\alpha_0$, if there exists an integer-valued random variable $N$ satisfying $\tau_N=T$, $P$-a.s. and $N\in L^{\infty}(\mathcal{F}_T)$.
\end{definition}
We denote by
\begin{itemize}
    \item $\mathcal{M}$: the set of all admissible switching strategies;
    \item $\mathcal{M}^i$: the subset of $\mathcal{M}$ consisting of admissible switching strategies starting from mode $i$;
    \item $\mathcal{M}_t$: the set of all admissible strategies starting at time $t$ (or equivalently $\tau_0=t$);
    \item $\mathcal{M}^i_t$: the subset of $\mathcal{M}_t$ consisting of admissible switching strategies starting at time $t$ from mode $i$.
\end{itemize}
For any $a\in\mathcal{M}_t$, we introduce the process $A^{a}$ as follows
\begin{equation*}
    A^{a}(s)=\sum_{j=1}^{N-1}k_{\alpha_{j-1},\alpha_j}\mathbf{1}_{[\tau_j,T]}(s),~~s\in[t,T]
\end{equation*}
which implies that $A^{a}$ is an adapted non-decreasing c\`{a}dl\`{a}g process and $A^{a}(T)\in L^{\infty}(\mathcal{F}_T)$ due to the fact that $N\in L^{\infty}(\mathcal{F}_T)$.

Since $A^a(t)$ is bounded for all $t\in[0,T]$, in view of (A1), it follows from \cite[Theorem 2.2.]{BE1} that the following BSDE 
\begin{equation}\label{eq:switchedtf_BSDE}
    \tilde{U}(s)=\xi_{a(T)}+A^a(T)+\int_s^Tg_{a(r)}(r,\tilde{U}(r)-A^a(r),\tilde{V}(r))dr-\int_s^T\tilde{V}(r)dW(r),~~~s\in[t,T].
\end{equation} has a unique solution $(\tilde{U},\tilde{V})$ such that $(\tilde{U},\tilde{Z}\cdot W)\in\mathcal{S}^{\infty}(\mathbb{R}^n)\times BMO$. Therefore, by a change of variable, 
\begin{equation*}
    U^a(s)=\tilde{U}(s)-A^a(s),~~V^a(s)=\tilde{V}(s),~~~s\in[0,T],
\end{equation*}
the following switched BSDE
\begin{equation}\label{eq:switched_BSDE}
    U^a(s)=\xi_{a(T)}+A^a(T)-A^a(s)+\int_s^Tg_{a(r)}(r,U^a(r),V^a(r))dr-\int_s^TV^a(r)dW(r),~~~s\in[t,T]
\end{equation}
admits a unique soltion $(U^a,V^a)$ such that $(U^a,V^a\cdot W)\in\mathcal{S}^{\infty}(\mathbb{R}^n)\times BMO$. 

For any $i=1,2,\ldots,n$, the optimal switching problem with initial mode $i$ is to minimize $U^a(t)$ subject to $a\in\mathcal{M}^i_t$. 

For a solution $(Y,Z,K)$ of RBSDE \eqref{BSDE} satisfying $(Y,Z\cdot W,K)\in\mathcal{S}^{\infty}_c(\mathbb{R}^n)\times BMO \times\mathcal{A}^{\infty}_c(\mathbb{R}^n)$, by setting $\tau^*_0=t$, $\alpha^*_0=i$, we define $(\tau^*_j,\alpha^*_j)$ for $j=1,2,\ldots$, in an inductive way as follows
\begin{equation*}
    \tau^*_j:=\inf\left\{t\geq \tau^*_{j-1}:Y_{\alpha^*_{j-1}}(t)=\min_{l\neq \alpha^*_{j-1}}\left\{Y_l(s)+k_{\alpha^*_{j-1},l}\right\}\right\}\wedge T,
\end{equation*}
and $\alpha^*_j$ is $\mathcal{F}_{\tau^*_j}$-measurable random variable such that 
\begin{equation*}
    Y_{\alpha^*_{j-1}}(\tau^*_j)=Y_{\alpha^*_j}(\tau^*_j)+k_{\alpha^*_{j-1},\alpha^*_j}.
\end{equation*}
We have the following representation result for $Y$, which immediately implies the uniqueness part of Theorem \ref{thm:main}.
\begin{theorem}\label{thm:uniqueness}
    Assuming that (A1), (A2) and (A4) hold, then
    \begin{itemize}
        \item[(i)] For any $a\in\mathcal{M}^i_t$, it holds that
                \begin{equation*}
                    \tilde{Y}_i(t)\leq U^{a}(t),~~~P\text{-a.s.}
                \end{equation*} 
        \item[(ii)] There exists an integer-valued random variable $N^*$ such that $\tau^*_{N^*}=T$, $P$-a.s. and $N^*\in L^{\infty}(\mathcal{F}_T)$. Morever, the following switching strategy
                \begin{equation*}
                    a^*(s)=i\mathbf{1}_{\{t\}}(s)+\sum_{j=1}^{N^*}\alpha^*_{j-1}\mathbf{1}_{(\tau^*_{j-1},\tau^*_{j}]}(s)
                \end{equation*}
                is admissible, i.e., $a^*\in\mathcal{M}^i_t$, and it holds that
                \begin{equation*}
                    Y_i(t)=U^{a^*}(t).
                \end{equation*}
        \item[(iii)] $Y$ satisfies the following representation
                \begin{equation*}
                    Y_i(t)=\essinf_{a\in\mathcal{M}^i_t}U^a(t),~~i=1,2,\ldots,n,~t\in[0,T].
                \end{equation*}
            Moreover, RBSDE \eqref{BSDE} has a unique solution $(Y,Z,K)$ such that $(Y,Z\cdot W,K)\in\mathcal{S}^{\infty}_c(\mathbb{R}^n)\times BMO \times\mathcal{A}^{\infty}_c(\mathbb{R}^n)$.
    \end{itemize}
\end{theorem}
\begin{proof}
    Without loss of generality, it is sufficient to prove (i) and (ii) for the case of $t=0$. Otherwise, we could consider the admissible switching strategies starting at time $t$.

    (i) We first introduce a triplet of processes $(Y^a,Z^a,K^a)$:
    \begin{equation*}
        Y^a(s)=\sum_{i=1}^NY_{\alpha_{i-1}}(s)\mathbf{1}_{[\tau_{i-1},\tau_i)}(s)+\xi_{a(T)}\mathbf{1}_{\{T\}}(s),
    \end{equation*}
    \begin{equation*}
        Z^a(s)=\sum_{i=1}^NZ_{\alpha_{i-1}}(s)\mathbf{1}_{[\tau_{i-1},\tau_i)}(s),
    \end{equation*}
    \begin{equation*}
        K^a(s)=\sum_{i=1}^N\int_{\tau_{i-1}\wedge s}^{\tau_i\wedge s}dK_{\alpha_{i-1}}(s).
    \end{equation*}
    Since $(Y,Z\cdot W,K)\in\mathcal{S}^{\infty}_c(\mathbb{R}^n)\times BMO \times\mathcal{A}^{\infty}_c(\mathbb{R}^n)$, one could easily check that $(Y^a,Z^a\cdot W,K^a)\in\mathcal{S}^{\infty}(\mathbb{R}^n)\times BMO \times\mathcal{A}^{\infty}(\mathbb{R}^n)$.
    Noting that $Y^a$ is a c\`{a}dl\`{a}g process which has jump $Y_{\alpha_i}(\theta_i)-Y_{\alpha_{i-1}}(\theta_i)$ at $\theta_i$ for $i=1,2,\ldots,N-1$, we get
    \begin{align*}
        &Y^a(s)-Y^a(0)\\
        &=-\sum_{i=1}^N\int_{\tau_{i-1}\wedge s}^{\tau_i\wedge s}g_{\alpha_{i-1}}(r,Y_{\alpha_{i-1}}(r),Z_{\alpha_{i-1}}(r))dr+\sum_{i=1}^N\int_{\tau_{i-1}\wedge s}^{\tau_i\wedge s}Z_{\alpha_{i-1}}(r)dW(r)+\sum_{i=1}^N\int_{\tau_{i-1}\wedge s}^{\tau_i\wedge s}dK_{\alpha_{i-1}}(r)\\
        &\qquad +\sum_{i=1}^{N-1}\left(Y_{\alpha_i}(\tau_i)-Y_{\alpha_{i-1}}(\tau_i)\right)\mathbf{1}_{[\tau_i,T]}(s)\\
        &=-\int_{0}^{s}g_{a(r)}(r,Y^a(r),Z^a(r))dr+\int_{0}^{s}Z^a(r)dW(r)+\int_{0}^{s}dK^a(r)+\tilde{A}^a(s)-A^a(s)
    \end{align*}
    where 
    \begin{equation*}
        \tilde{A}^a(s)=\sum_{i=1}^{N-1}\left(Y_{\alpha_i}(\tau_i)+k_{\alpha_{i-1},\alpha_i}-Y_{\alpha_{i-1}}(\tau_i)\right)\mathbf{1}_{[\tau_i,T]}(s).
    \end{equation*}
    Since $Y\in\mathcal{S}^{\infty}(\mathbb{R}^n)$, $N\in L^{\infty}(\mathcal{F}_T)$ and $Y(t)\in\overline{Q}$ for any $t\in[0,T]$, we deduce that $A^a$ is an adapted non-decreasing c\`{a}dl\`{a}g process satisfying $A^a(T)\in L^{\infty}(\mathcal{F}_T)$. Therefre, we obtain that $(Y^a,Z^a)$ is a solution of the following BSDE
    \begin{align*}
        Y^a(s)&=\xi_{a(T)}+A^a(T)-A^a(s)-\left(\left(K^a(T)+\tilde{A}^a(T)\right)-\left(K^a(s)+\tilde{A}^a(s)\right)\right)\\
        &\quad +\int_s^Tg_{a(r)}(r,Y^a(r),Z^a(r))dr-\int_s^TZ^a(r)dW(r),~~~s\in[0,T].
    \end{align*}
    Since both $K^a$ and $\tilde{A}^a$ are adapted non-decreasing c\`{a}dl\`{a}g processes, it follows from \cite[Theorem 2]{Te} that 
    \begin{equation*}
        Y^a(0)\leq U^a(0),
    \end{equation*}
    from which, we obtain 
    \begin{equation*}
        Y_i(0)\leq U^a(0).
    \end{equation*}
    (ii) We first introduce the following (closed) subsets of $\overline{Q}$: for $i\neq j$,
    \begin{equation*}
        B_{i,j}:=\left\{y\in\mathbb{R}^n:~y_i=y_j+k_{i,j}\right\}\cap\overline{Q}.
    \end{equation*}
    Since if there exists $y\in\overline{Q}$ satisfying for $i\neq j$ and $j\neq l$  
    \begin{equation*}
        y_i=y_j+k_{i,j}~\text{and}~y_j=y_l+k_{jl}
    \end{equation*}
    it holds in view of (A4) that 
    \begin{equation*}
        y_i=y_l+k_{i,j}+k_{j,l}>y_l+k_{i,l}
    \end{equation*}
    which contradicts the fact that $y\in\overline{Q}$, we deduce that for $i\neq j$ and $j\neq l$, $B_{i,j}\cap B_{j,l}=\emptyset$. Therefore, the distance between $B_{i,j}$ and $B_{j,l}$ is strictly poisitve, i.e.,
    \begin{equation*}
        dist\left(B_{i,j},B_{j,l}\right)>0.
    \end{equation*}
    By setting
    \begin{equation*}
        c:=\min_{i\neq j,j\neq l}dist\left(B_{i,j},B_{j,l}\right)>0,
    \end{equation*}
    for $0\leq \tau^*_1<\tau^*_2<T$, it follows from the definition of $(\tau^*_1,\alpha^*_1)$ and $(\tau^*_2,\alpha^*_2)$ that 
    \begin{equation*}
        Y(\tau^*_1)\in B_{i,\alpha^*_1}~~\text{and}~~Y(\tau^*_2)\in B_{\alpha^*_1,\alpha^*_2},
    \end{equation*}
    which implies that 
    \begin{equation*}
        |Y(\tau^*_2)-Y(\tau^*_1)|\geq c.
    \end{equation*}
    Similarly, if $\tau^*_1<\tau^*_2<\cdots<\tau^*_{j-1}<\tau^*_j<T$, it holds that 
    \begin{equation}\label{eq:dist}
        |Y(\tau^*_j)-Y(\tau^*_j-1)|\geq c.
    \end{equation}
    On the other hand, since $(Y,Z,K)$ is a solution of RBSDE \eqref{BSDE} satisfying $(Y,Z\cdot W,K)\in\mathcal{S}^{\infty}_c(\mathbb{R}^n)\times BMO \times\mathcal{A}^{\infty}_c(\mathbb{R}^n)$, one could easily check that 
    \begin{equation}\label{eq:bound}
        \left\|\sum_{j=1}^{\infty}|Y(\tau^*_j)-Y(\tau^*_{j-1})|\right\|_{\infty}<\infty.
    \end{equation}
    Therefore, defining 
    \begin{equation*}
        N^*=\inf\left\{j|\tau^*_j=T\right\},
    \end{equation*}
    combining \eqref{eq:dist} and \eqref{eq:bound}, it holds that
    \begin{equation*}
      N^*\in L^\infty(\mathcal{F}_T),~~\tau^*_{N^*}=T,~~P\text{-a.s.}  
    \end{equation*}
    Thus, we deduce that $A^{a^*}\in L^{\infty}(\mathcal{F}_T)$ and $a^*$ is admissible. Hence, we could introduce the processes $A^{a^*},Y^{a^*},Z^{a^*},K^{a^*}$ and $\tilde{A}^{a^*}$ similarly as above. On the other hand, from the choice of $a^*$, we obtain
    \begin{equation*}
        K^{a^*}+\tilde{A}^{a^*}=0.
    \end{equation*}
    Hence, we deduce that 
    \begin{equation*}
        Y^{a^*}(s)=\xi_{a^*(T)}+A^{a^*}(T)-A^{a^*}(s)+\int_t^Tg_{a^*(r)}(r,Y^{a^*}(r),Z^{a^*}(r))dr-\int_s^TZ^{a^*}(r)dW(r),~~s\in[0,T]
    \end{equation*}
    which implies
    \begin{equation*}
        Y^{a^*}(0)=U^{a^*}(0)
    \end{equation*}
    or equivalently
    \begin{equation*}
        Y_i(0)=U^{a^*}(0).
    \end{equation*}
    (iii) The representation for $Y$ follows immediately from both assertions (i) and (ii), which implies the uniqueness of the first component of the solution of RBSDE \eqref{BSDE}. As the consequence of a direct computation, the uniqueness of the other two components of the solution of RBSDE \eqref{BSDE} follows.
\end{proof}
\subsection{Proof of Theorem \ref{thm:main-ex}}
\label{extension}

In veiw of (A5), for any $y\in\mathcal{S}^{\infty}_{c}(\mathbb{R}^n)$, it follows from Theorem \ref{thm:main} that the following RBSDE
\begin{equation}\label{BSDE2}
    \begin{cases}
    &Y_i(t)=\xi_i+\int_t^Tg_i(s,y(s),Z_i(s))ds-\int_t^TdK_i(s)-\int_t^TZ_i(s)dW(s),~~t\in[0,T],\\
    &Y_i(t)\leq \min\limits_{j\neq i}\{Y_j(t)+k_{i,j}\},\\
    &\int_0^T\left(Y_i(s)-\min\limits_{j\neq i}\left\{Y_j(s)+k_{i,j}\right\}\right)dK_i(s)=0,~~i=1,2,\ldots,n,
    \end{cases}
    \end{equation}
    adimits a unique solution $(Y,Z,K)$ satisfying $(Y,Z\cdot W,K)\in\mathcal{S}^{\infty}_c(\mathbb{R}^n)\times BMO \times\mathcal{A}^{\infty}_c(\mathbb{R}^n)$. We now introduce a functional $\Phi$ mapping from $\mathcal{S}^{\infty}_{c}(\mathbb{R}^n)$ to $\mathcal{S}^{\infty}_c(\mathbb{R}^n)$ by $\Phi(y):=Y$, where $Y$ is the first component of the unique solution of RBSDE \eqref{BSDE2} associated to $y$.

    Meanwhile, for constant $\beta>0$, we introduce the following equivalent norm $\|\cdot\|_{\beta,\mathcal{S}^{\infty}}$ on $\mathcal{S}^{\infty}_c(\mathbb{R}^n)$ by 
    \begin{equation*}
        \|Y\|_{\beta,\mathcal{S}^{\infty}}:=\bigg\| \sup_{0\leq t\leq T} e^{\beta t} |Y(t)| \bigg\|_{\infty}.
    \end{equation*}
    Now for $y,\bar{y}\in\mathcal{S}^{\infty}_c(\mathbb(R)^n)$ and any $i=1,2,\ldots,n$, we define function $\tilde{g}_i$ as follows
    \begin{equation*}
        \tilde{g}_i(t,z):=g_i(t,y(t),z)\wedge g_i(t,\bar{y}(t),z),~~~t\in[0,T], z\in\mathbb{R}^d.
    \end{equation*}
    Since (A5) holds, again it follows from Theorem \ref{thm:main} that the following RBSDE 
    \begin{equation}\label{BSDE3}
        \begin{cases}
        &\tilde{Y}_i(t)=\xi_i+\int_t^T\tilde{g}_i(s,\tilde{Z}_i(s))ds-\int_t^Td\tilde{K}_i(s)-\int_t^T\tilde{Z}_i(s)dW(s),~~t\in[0,T],\\
        &\tilde{Y}_i(t)\leq \min\limits_{j\neq i}\{\tilde{Y}_j(t)+k_{i,j}\},\\
        &\int_0^T\left(\tilde{Y}_i(s)-\min\limits_{j\neq i}\left\{\tilde{Y}_j(s)+k_{i,j}\right\}\right)d\tilde{K}_i(s)=0,~~i=1,2,\ldots,n,
        \end{cases}
        \end{equation}
        adimits a unique solution $(\tilde{Y},\tilde{Z},\tilde{K})$ satisfying $(\tilde{Y},\tilde{Z}\cdot W,\tilde{K})\in\mathcal{S}^{\infty}_c(\mathbb{R}^n)\times BMO \times\mathcal{A}^{\infty}_c(\mathbb{R}^n)$.
    Denoting
    \begin{equation*}
        Y=\Phi(y)~~~\text{and}~~~ \bar{Y}=\Phi(\bar{y}),
    \end{equation*}  
    it follows from Theorem \ref{thm:uniqueness} that for $t\in[0,T]$ and $i=1,2,\ldots,n$,
    \begin{equation}\label{eq:rep-ex}
        Y_i(t)=\essinf_{a\in\mathcal{M}^i_t}U^a(t),~~\bar{Y}_i(t)=\essinf_{a\in\mathcal{M}^i_t}\bar{U}^a(t),~~\text{and}~~\tilde{Y}_i(t)=\essinf_{a\in\mathcal{M}^i_t}\tilde{U}^a(t),
    \end{equation} 
    where 
    \begin{align*}
        &U^a(s)=\xi_{a(T)}+A^a(T)-A^a(s)+\int_s^Tg_{a(r)}(r,y(r),V^a(r))dr-\int_s^TV^a(r)dW(r),~~~s\in[t,T]\\
        &\bar{U}^a(s)=\xi_{a(T)}+A^a(T)-A^a(s)+\int_s^Tg_{a(r)}(r,\bar{y}(r),\bar{V}^a(r))dr-\int_s^T\bar{V}^a(r)dW(r),~~~s\in[t,T]\\
        &\tilde{U}^a(s)=\xi_{a(T)}+A^a(T)-A^a(s)+\int_s^T\tilde{g}_{a(r)}(r,\tilde{V}^a(r))dr-\int_s^T\tilde{V}^a(r)dW(r),~~~s\in[t,T].
    \end{align*}
    Moreover, for $t\in[0,T]$ and $i=1,2,\ldots,n$, there exists $\tilde{a}^*\in\mathcal{M}^i_t$ such that
    \begin{equation}\label{eq:rep1-ex}
        \tilde{Y}_i(t)=\tilde{U}^{\tilde{a}^*}(t).
    \end{equation}
    Thus, using Theorem 2 of \cite{Te}, we get
    \begin{equation*}
        \tilde{U}^a(s)\leq U^a(s),~~\text{and}~~\tilde{U}^a(s)\leq \bar{U}^a(s),~~~\text{for any }~a\in\mathcal{M}^i_t,~s\in[t,T].
    \end{equation*}
    Therefore, in view of \eqref{eq:rep-ex} and \eqref{eq:rep1-ex} and noting that $\tilde{a}^*\in\mathcal{M}^i_t$, it holds that 
    \begin{equation*}
        \tilde{U}^{\tilde{a}^*}(t)\leq Y_i(t)\leq U^{\tilde{a}^*}(t),~~\text{and}~~\tilde{U}^{\tilde{a}^*}(t)\leq \bar{Y}_i(t)\leq \bar{U}^{\tilde{a}^*}(t),~~~\text{for any }~i=1,2,\ldots,n,~t\in[0,T],
    \end{equation*}
    from which, we deduce that 
    \begin{equation}\label{eq:estimate-ex}
        |Y_i(t)-\bar{Y}_i(t)|\leq |U^{\tilde{a}^*}(t)-\tilde{U}^{\tilde{a}^*}(t)|+|\bar{U}^{\tilde{a}^*}(t)-\tilde{U}^{\tilde{a}^*}(t)|.
    \end{equation}
    Applying It\^{o}'s formula to $e^{\beta t}|U^{\tilde{a}^*}(t)-\tilde{U}^{\tilde{a}^*}(t)|^2$ and using (A5) and inequality $|x\wedge y-y|\leq |x-y|$, we obtain
    \begin{align*}
        &e^{\beta t}|U^{\tilde{a}^*}(t)-\tilde{U}^{\tilde{a}^*}(t)|^2+\int_t^Te^{\beta s}|V^{\tilde{a}^*}(s)-\tilde{V}^{\tilde{a}^*}(s)|^2ds\\
        &=-\beta\int_t^Te^{\beta s}|U^{\tilde{a}^*}(s)-\tilde{U}^{\tilde{a}^*}(s)|^2ds-2\int_t^Te^{\beta s}\left(U^{\tilde{a}^*}(s)-\tilde{U}^{\tilde{a}^*}(s)\right)\left(V^{\tilde{a}^*}(s)-\tilde{V}^{\tilde{a}^*}(s)\right)dW(s)\\
        &+2\int_t^Te^{\beta s}\left(U^{\tilde{a}^*}(s)-\tilde{U}^{\tilde{a}^*}(s)\right)\left(g_{\tilde{a}^*(s)}(s,y(s),V^{\tilde{a}^*}(s))-\tilde{g}_{\tilde{a}^*(s)}(s,\tilde{V}^{\tilde{a}^*}(s))\right)ds\\
        &\leq-\beta\int_t^Te^{\beta s}|U^{\tilde{a}^*}(s)-\tilde{U}^{\tilde{a}^*}(s)|^2ds-2\int_t^Te^{\beta s}\left(U^{\tilde{a}^*}(s)-\tilde{U}^{\tilde{a}^*}(s)\right)\left(V^{\tilde{a}^*}(s)-\tilde{V}^{\tilde{a}^*}(s)\right)d\tilde{W}^*(s)\\
        &+2\gamma\int_t^Te^{\beta s}|U^{\tilde{a}^*}(s)-\tilde{U}^{\tilde{a}^*}(s)||y(s)-\bar{y}(s)|ds
    \end{align*}
    where $\tilde{W}^*(s)=W(s)-\int_0^s\tilde{\beta}^*(r)dr$ is a Brownian motion under an equivalent probability measure $\frac{d\tilde{P}^*}{dP}=\mathcal{E}_T\left(\tilde{\beta}^*\cdot W\right)$, and the process $\tilde{\beta}^{*}$ satisfies $|\tilde{\beta}^*|\leq\gamma(1+|V^{\tilde{a}^*}|+|\tilde{V}^{\tilde{a}^*}|)$. Taking conditional expectation with respect to  $\mathcal{F}_t$ and $\tilde{E}^*$ and using Young's inequality, we obtain
    \begin{align*}
        e^{\beta t}|U^{\tilde{a}^*}(t)-\tilde{U}^{\tilde{a}^*}(t)|^2&\leq -\beta\tilde{E}^*\left[\int_t^Te^{\beta s}|U^{\tilde{a}^*}(s)-\tilde{U}^{\tilde{a}^*}(s)|^2ds\bigg|\mathcal{F}_t\right]\\
        &\quad +2\gamma \tilde{E}^*\left[\int_t^Te^{\beta s}|U^{\tilde{a}^*}(s)-\tilde{U}^{\tilde{a}^*}(s)||y(s)-\bar{y}(s)|ds\bigg|\mathcal{F}_t\right]\\
        &\leq \frac{\gamma^2}{\beta} \tilde{E}^*\left[\int_t^Te^{\beta s}|y(s)-\bar{y}(s)|^2ds\bigg|\mathcal{F}_t\right]\\
        &\leq \frac{\gamma^2T}{\beta}\|y-\bar{y}\|^2_{\beta,\mathcal{S}^{\infty}}.
    \end{align*}
    Similarly, it holds that 
    \begin{equation*}
        e^{\beta t}|U^{\tilde{a}^*}(t)-\tilde{U}^{\tilde{a}^*}(t)|^2\leq \frac{\gamma^2T}{\beta}\|y-\bar{y}\|^2_{\beta,\mathcal{S}^{\infty}}.
    \end{equation*}
    Recalling \eqref{eq:estimate-ex}, we deduce that 
    \begin{equation*}
        e^{\beta t}|Y_i(t)-\bar{Y}_i(t)|^2\leq \frac{2\gamma^2T}{\beta}\|y-\bar{y}\|^2_{\beta,\mathcal{S}^{\infty}}.
    \end{equation*}
    Thus, it holds that 
    \begin{equation*}
        \|Y-\bar{Y}\|^2_{\beta,\mathcal{S}^{\infty}}\leq \frac{2\gamma^2Tn}{\beta}\|y-\bar{y}\|^2_{\beta,\mathcal{S}^{\infty}}.
    \end{equation*}
    Hence, by choosing $\beta>2\gamma^2 Tn$, we obtain that $\Phi$ is a contraction mapping on $\mathcal{S}^{\infty}_c(\mathbb{R}^n)$ which concludes the proof.
\section{Risk-sensitive optimal switching of functional SDEs}
In this section, we consider a risk-sensitive optimal switching of functional SDEs. We consider functions $\sigma:[0,T]\times\mathcal{C}_T(\mathbb{R}^d)\rightarrow\mathbb{R}^{d\times d}$, $b:[0,T]\times\mathcal{C}_T(\mathbb{R}^d)\rightarrow \mathbb{R}^{n\times d}$ and $l:[0,T]\times\mathcal{C}_T(\mathbb{R}^d)\rightarrow \mathbb{R}^{n\times d}$. We will make the following assumptions.
\begin{itemize}
    \item[(H1)] For any continuous $\mathcal{P}$-measurable process $x=(x_t)_{0\leq t\leq T}$, processes $(\sigma(t,x))_{0\leq t\leq T}$, $(b(t,x))_{0\leq t\leq T}$ and $(l(t,x))_{0\leq t\leq T}$ are $\mathcal{P}$-measurable.
    \item[(H2)] There exists a constant $\beta>0$ such that for $t\in[0,T]$ and $x,\bar{x}\in\mathcal{C}_T(\mathbb{R}^d)$,
        \begin{equation*}
            |b(t,x)-b(t,\bar{x})|+|\sigma(t,x)-\sigma(t,\bar{x})|+|l(t,x)-l(t,\bar{x})|\leq \beta \sup_{0\leq s\leq t}|x_s|.
        \end{equation*} 
    \item[(H3)] $b$ and $l$ are bounded, and $\int_0^T|\sigma(t,0)|^2dt<\infty$.
\end{itemize}

Let $x_0$ be a fixed point in $\mathbb{R}^d$. Under the above assumptions, the following functional SDE:
\begin{equation}\label{eq:functionalSDE}
    X(t)=x_0+\int_0^t\sigma(s,X)dW(s),~~t\in[0,T]
\end{equation}
has a unique strong solution $X$ in $\mathcal{S}^2(\mathbb{R}^d)$.
 
For each $a\in\mathcal{M}^i$, we consider the switched equation
\begin{equation}\label{eq:switchSDE}
    X^a(t)=x_0+\int_0^t\sigma(s,X^a)\left(dW(s)+b_{a(s)}(s,X^a)ds\right),~t\in[0,T].
\end{equation}
In view of (H3), the process
\begin{equation*}
    W^{a}(t)=W(t)-\int_0^tb_{a(s)}(s,X)ds,~~~t\in[0,T]
\end{equation*}
is a Brownian motion under an equivalent probability measure 
\begin{equation*}
    \frac{dP^{a}}{dP}=\mathcal{E}_T\left(b_{a(\cdot)}(\cdot,X)\cdot W\right).
\end{equation*}
Since $X$ is the strong solution of SDE \eqref{eq:functionalSDE}, we dedcue that the triplet $(P^{a},W^{a},X)$ is a weak solution of SDE \eqref{eq:switchSDE}.

We consider the cost functional
\begin{equation*}
    J(a)=E^a\left[\exp\left(\int_0^Tl_{a(s)}(s,X)ds+\sum_{i=1}^{N-1}k_{\alpha_{i-1},\alpha_{i}}+\xi_{a(T)}\right)\right].
\end{equation*}
The switching problem is to minimize the cost $J(a)$ over $a\in\mathcal{M}^i$, subject to the state equation.
For any $i=1,2,\ldots,n$, we define $\varphi_i$ as follows: $\forall (t,x,z)\in[0,T]\times \mathcal{C}_T(\mathbb{R}^d)\times\mathbb{R}^d$,
\begin{equation*}
    \varphi_i(t,x,z):=l_i(t,x)+\langle z,b_i(t,x)\rangle+\frac{1}{2}z^2.
\end{equation*}
Under the assumptions (A2), (A4), (H1) and (H3), it follows from Theorem \ref{thm:main} that the following RBSDE
\begin{equation}\label{eq:riskBSDE}
    \begin{cases}
    &Y_i(t)=\xi_i+\int_t^T\varphi_i(s,X,Z_i(s))ds-\int_t^TdK_i(s)-\int_t^TZ_i(s)dW(s),~~t\in[0,T],\\
    &Y_i(t)\leq \min\limits_{j\neq i}\{Y_j(t)+k_{i,j}\},\\
    &\int_0^T\left(Y_i(s)-\min\limits_{j\neq i}\left\{Y_j(s)+k_{i,j}\right\}\right)dK_i(s)=0, ~~i=1,2,\ldots,n,
    \end{cases}
    \end{equation}
    admits a unique solution $(Y,Z,K)$ such that $(Y,Z\cdot W,K)\in\mathcal{S}^{\infty}_{c}(\mathbb{R}^n)\times BMO\times\mathcal{S}^{\infty}_{c}(\mathbb{R}^n)$.
\begin{theorem}
    Suppose that (A2), (A4), (H1), (H2), (H3) hold and $1\leq i\leq n$. Then
    \begin{itemize}
        \item[(i)] For each $a\in\mathcal{M}^i$ and the associated weak solution $(P^a,W^a,X)$ of SDE \eqref{eq:switchSDE}, we have
        \begin{equation*}
            J(a)\geq e^{Y_i(0)}.
        \end{equation*} 
        \item[(ii)] There exists an optimal switching strategy $a^*\in\mathcal{M}^i$ and the associated weak solution $(P^{a^*},W^{a^*},X)$ of SDE \eqref{eq:switchSDE} such that 
        \begin{equation*}
            J(a^*)=e^{Y_i(0)}.
        \end{equation*} 
    \end{itemize}
\end{theorem}
\begin{proof}
    (i) For any $a\in\mathcal{M}^i$, in view of (H3), the process
    \begin{equation*}
        \tilde{W}(t)=W^{a}(t)-\int_0^tb_{a(s)}(s,X)ds,~~~t\in[0,T]
    \end{equation*}
    is a Brownian motion under an equivalent probability measure 
    \begin{equation*}
        \frac{d\tilde{P}}{dP^a}=\mathcal{E}_T\left(b_{a(\cdot)}(\cdot,X)\cdot W^{a}\right).
    \end{equation*}
    Therefore, it follows from Theorem \ref{thm:main} that the following RBSDE
     \begin{equation}\label{eq:riskBSDE1}
        \begin{cases}
        &\tilde{Y}_i(t)=\xi_i+\int_t^T\varphi_i(s,X,\tilde{Z}_i(s))ds-\int_t^Td\tilde{K}_i(s)-\int_t^T\tilde{Z}_i(s)d\tilde{W}(s),~~t\in[0,T],\\
        &\tilde{Y}_i(t)\leq \min\limits_{j\neq i}\{\tilde{Y}_j(t)+k_{i,j}\},\\
        &\int_0^T\left(\tilde{Y}_i(s)-\min\limits_{j\neq i}\left\{\tilde{Y}_j(s)+k_{i,j}\right\}\right)d\tilde{K}_i(s)=0, ~~i=1,2,\ldots,n,
        \end{cases}
        \end{equation}
        admits a unique solution $(\tilde{Y},\tilde{Z},\tilde{K})$ such that $(\tilde{Y},\tilde{Z}\cdot \tilde{W},\tilde{K})\in\mathcal{S}^{\infty}_{c}(\mathbb{R}^n;\tilde{P})\times BMO(\tilde{P})\times\mathcal{S}^{\infty}_{c}(\mathbb{R}^n;\tilde{P})$. Thus by a classical argument of Yamada-Wanatane, for RBSDE \eqref{eq:riskBSDE}, the pathwise uniqueness implies the uniqueness in the sense of probability law \footnote{Indeed, by a standard exponential transform, it follows from the argument for Lipschitz BSDE, see \cite[Remark 1.6]{De}}. Hence, we have
    \begin{equation*}
        Y_i(0)=\tilde{Y}_i(0),~~i=1,2,\ldots,n.
    \end{equation*}
    We now define a quadruple of processes $(\tilde{Y}^a,\tilde{Z}^a,\tilde{K}^a,\tilde{A}^a)$ as follows:
    \begin{equation*}
        \tilde{Y}^a(s)=\sum_{i=1}^N\tilde{Y}_{\alpha_{i-1}}(s)\mathbf{1}_{[\tau_{i-1},\tau_i)}(s)+\xi_{a(T)}\mathbf{1}_{\{T\}}(s),
    \end{equation*}
    \begin{equation*}
        \tilde{Z}^a(s)=\sum_{i=1}^N\tilde{Z}_{\alpha_{i-1}}(s)\mathbf{1}_{[\tau_{i-1},\tau_i)}(s),
    \end{equation*}
    \begin{equation*}
        \tilde{K}^a(s)=\sum_{i=1}^N\int_{\tau_{i-1}\wedge s}^{\tau_i\wedge s}d\tilde{K}_{\alpha_{i-1}}(s),
    \end{equation*}
    \begin{equation*}
        \tilde{A}^a(s)=\sum_{i=1}^{N-1}\left(\tilde{Y}_{\alpha_i}(\tau_i)+k_{\alpha_{i-1},\alpha_i}-\tilde{Y}_{\alpha_{i-1}}(\tau_i)\right)\mathbf{1}_{[\tau_i,T]}(s).
    \end{equation*}
    Recalling the cost process $A^a$, we obtain in a similar way as in the proof of Theorem \ref{thm:uniqueness} that $(\tilde{Y}^a,\tilde{Z}^a)$ satisfies $(\tilde{Y}^a,\tilde{Z}^a\cdot\tilde{W})\in\mathcal{S}^{\infty}(\mathbb{R};\tilde{P})\times BMO(\tilde{P})$ and the following BSDE
    \begin{align*}
        \tilde{Y}^a(t)&=\xi_{a(T)}+A^a(T)-A^a(t)-\left(\left(\tilde{K}^a(T)+\tilde{A}^a(T)\right)-\left(\tilde{K}^a(t)+\tilde{A}^a(t)\right)\right)\\
        &+\int_t^T\varphi_{a(s)}(s,X,\tilde{Z}^a(s))ds-\int_t^T\tilde{Z}^a(s)d\tilde{W}(s),~~t\in[0,T].
    \end{align*}
    Since $\tilde{K}^a$ and $\tilde{A}^a$ are non-decreasing, it follows from a standard exponential transform that
    \begin{align*}
        e^{\tilde{Y}^a(0)}&\leq E^a\left[\exp\left(\xi_{a(T)}+A^a(T)+\int_0^Tl_{a(t)}(t,X)dt\right)\right].
    \end{align*}
    On the other hand, it holds from the definition of $\tilde{Y}^a$ that
    \begin{equation*}
        \tilde{Y}_i(0)=\tilde{Y}^a(0).
    \end{equation*}
    Therefore, we deduce that
    \begin{equation*}
        J(a(\cdot))\geq e^{Y_i(0)}.
    \end{equation*}
    (ii) Let $X$ be the unique strong solution in $\mathcal{S}^2(\mathbb{R}^d)$ of SDE \eqref{eq:functionalSDE} and $(Y,Z,K)$ be the unique solution of RBSDE \eqref{eq:riskBSDE} satisfying $(Y,Z\cdot W,K)\in\mathcal{S}^{\infty}_c(\mathbb{R}^n)\times BMO \times\mathcal{A}^{\infty}_c(\mathbb{R}^n)$. By setting $\tau^*_0=0$, $\alpha^*_0=i$, we define $(\tau^*_j,\alpha^*_j)$ for $j=1,2,\ldots$, in an inductive way as follows
    \begin{equation*}
        \tau^*_j:=\inf\left\{t\geq \tau^*_{j-1}:Y_{\alpha^*_{j-1}}(t)=\min_{l\neq \alpha^*_{j-1}}\left\{Y_l(s)+k_{\alpha^*_{j-1},l}\right\}\right\}\wedge T,
    \end{equation*}
    and $\alpha^*_j$ is $\mathcal{F}_{\tau^*_j}$-measurable random variable such that 
    \begin{equation*}
        Y_{\alpha^*_{j-1}}(\tau^*_j)=Y_{\alpha^*_j}(\tau^*_j)+k_{\alpha^*_{j-1},\alpha^*_j}.
    \end{equation*}
    Similarly as in the proof of Theorem \ref{thm:uniqueness}, there exists an integer-valued random variable $N^*$ such that
    \begin{equation*}
        N^*\in L^\infty(\mathcal{F}_T),~~\tau^*_{N^*}=T,~~P\text{-a.s.}  
      \end{equation*}
    Now we define the switching strategy $a^*$ as follows
    \begin{equation*}
        a^*(s)=i\mathbf{1}_{\{0\}}(s)+\sum_{j=1}^{N}\alpha^*_{j-1}\mathbf{1}_{(\tau^*_{j-1},\tau^*_{j}]}(s).
    \end{equation*}
    It is easy to verify that $A^{a^*}\in L^{\infty}(\mathcal{F}_T)$ and $a^*$ is admissible. On the other hand,  in view of (H3), the process
    \begin{equation*}
        W^{a^*}(t)=W(t)-\int_0^tb_{a^*(s)}(s,X)ds,~~~t\in[0,T]
    \end{equation*}
    is a Brownian motion under an equivalent probability measure 
    \begin{equation*}
        \frac{dP^{a^*}}{dP}=\mathcal{E}_T\left(b_{a^*(\cdot)}(\cdot,X)\cdot W\right).
    \end{equation*}
    Hence, the triplet $(P^{a^*},W^{a^*},X)$ is a weak solution of SDE \eqref{eq:switchSDE}, corresponding to $a^*$.

    Defining
    \begin{align*}
        A^{a^*}(s)&=\sum_{j=1}^{N^*-1}k_{\alpha^*_{j-1},\alpha^*_j}\mathbf{1}_{[\tau^*_j,T]}(s),\\
        Y^{a^*}(s)&=\sum_{i=1}^{N^*}Y_{\alpha^*_{i-1}}(s)\mathbf{1}_{[\tau^*_{i-1},\tau^*_i)}(s)+\xi_{a^*(T)}\mathbf{1}_{\{T\}}(s),\\
        Z^{a^*}(s)&=\sum_{i=1}^{N^*}Z_{\alpha^*_{i-1}}\mathbf{1}_{[\tau^*_{i-1},\tau^*_i)}(s),\\
        K^{a^*}(s)&=\sum_{i=1}^{N^*}\int_{\tau^*_{i-1}\wedge s}^{\tau^*_i\wedge s}dK_{\alpha^*_{i-1}}(r),\\
        \tilde{A}^{a^*}(s)&=\sum_{i=1}^{N^*-1}\left(Y_{\alpha_i}(\tau^*_i)+k_{\alpha^*_{i-1},\alpha^*_i}-Y_{\alpha^*_{i-1}}(\tau^*_i)\right)\mathbf{1}_{[\tau^*_i,T]}(s),
    \end{align*}
    we get
    \begin{equation*}
        (Y^{a^*},Z^{a^*}\cdot W)\in\mathcal{S}^{\infty}(\mathbb{R})\times BMO,~~K^{a^*}=0,~\tilde{A}^{a^*}=0.
    \end{equation*}
    Similarly as in the proof of Theorem \ref{thm:uniqueness}, we deduce that for any $t\in[0,T]$,
    \begin{align*}
        Y^{a^*}(t)&=\xi_{a^*(T)}+A^{a^*}(T)-A^{a^*}(t)-\left(\left(K^{a^*}(T)+\tilde{A}^{a^*}(T)\right)-\left(K^{a^*}(t)+\tilde{A}^{a^*}(t)\right)\right)\\
        &\quad+\int_t^T\varphi_{a^{*}(s)}(s,X,Z^{a^*}(s))ds-\int_t^TZ^{a^*}(s)dW(s)\\
        &=\xi_{a^*(T)}+A^{a^*}(T)-A^{a^*}(t)+\int_t^T\varphi_{a^{*}(s)}(s,X,Z^{a^*}(s))ds-\int_t^TZ^{a^*}(s)dW(s).
    \end{align*}
    It follows from a standard exponential transform that
    \begin{align*}
        e^{Y^{a^*}(0)}=E^{a^*}\left[\exp\left(\xi_{a^*(T)}+A^{a^*}(T)+\int_0^Tl_{a^*(t)}(t,X)dt\right)\right],
    \end{align*}
    which combining with the fact that $Y_i(0)=Y^{a^*}(0)$ implies that
    \begin{equation*}
        J(a^*)=e^{Y_i(0)}.
    \end{equation*}
\end{proof}

\end{document}